\documentclass[11pt,reqno]{amsart}
\usepackage{amssymb,amscd,amsbsy}
\usepackage{amssymb,amscd,amsbsy,mathrsfs}
\newcommand{\lb}{\linebreak}

\renewcommand{\a}{\alpha}

\newcommand{\e}{\varepsilon}
\newcommand{\vk}{\varkappa}
\newcommand{\z}{\zeta}

\renewcommand{\l}{\lambda}

\newcommand{\s}{\sigma}

\newcommand{\f}{\varphi}
\renewcommand{\o}{\omega}

\newcommand{\D}{\Delta}
\renewcommand{\L}{\Lambda}

\renewcommand{\O}{\Omega}

\newcommand{\E}{{\mathscr E}}
\newcommand{\cd}{{\mathscr D}}
\newcommand{\F}{{\mathscr F}}

\newcommand{\h}{{\mathscr H}}

\newcommand{\K}{{\mathscr K}}

\newcommand{\X}{{\mathscr X}}
\newcommand{\Y}{{\mathscr Y}}

\newcommand{\C}{{\Bbb C}}
\newcommand{\T}{{\Bbb T}}

\newcommand{\dd}{{\Bbb D}}
\newcommand{\R}{{\Bbb R}}
\newcommand{\Z}{{\Bbb Z}}

\newcommand{\0}{{\boldsymbol{0}}}

\newcommand{\bs}{\boldsymbol}

\newcommand{\bS}{{\boldsymbol S}}

\newcommand{\rf}[1]{(\ref{#1})}

\newcommand{\df}{\stackrel{\mathrm{def}}{=}}

\newcommand{\spn}{\operatorname{span}}
\newcommand{\supp}{\operatorname{supp}}
\newcommand{\clos}{\operatorname{clos}}

\newcommand{\const}{\operatorname{const}}

\newcommand{\eeq}{\end{equation}}
\newcommand{\beq}{\begin{equation}}
\newcommand{\bay}{\begin{eqnarray}}
\newcommand{\ba}{\begin{align*}}
\newcommand{\ea}{\end{align*}}
\newcommand{\ey}{\end{eqnarray}}
\newcommand{\bey}{\begin{eqnarray*}}
\newcommand{\eey}{\end{eqnarray*}}

\newcommand{\be}{\infty}

\newcommand{\bl}{\blacksquare}

\newcommand{\Range}{\operatorname{Range}}

\newcommand{\Pf}{{\bf Proof. }}
\newcommand{\im}{\operatorname{Im}}

\newtheorem{thm}{\hspace{\parindent}Theorem}[section]

\newtheorem{cor}[thm]{\hspace{\parindent}Corollary}
\newtheorem{lem}[thm]{\hspace{\parindent}Lemma}

\theoremstyle{remark}

\newtheorem*{rem*}{Remark}

\newcommand{\qm}{\quad\mbox{and}\quad}

\newcommand\Li{{\rm Lip}}
\newcommand\fM{\frak M}
\newcommand\cZ{\mathcal{Z}}
\newcommand\dg{\frak D}

\newcommand\ri{{\rm i}}

\begin{document}

\newcommand{\vse}{\vspace{.2in}}
\numberwithin{equation}{section}

\title{Functions of perturbed dissipative operators}
\author{A.B. Aleksandrov and V.V. Peller}
\thanks{The first author is partially supported by RFBR grant 08-01-00358-a;
the second author is partially supported by NSF grant DMS 1001844 and by ARC grant}

\newcommand{\mt}{{\mathcal T}}

\begin{abstract}
We generalize our results of \cite{AP2} and \cite{AP3} to the case of maximal dissipative operators. We obtain sharp conditions on a function analytic in the upper half-plane to be operator Lipschitz. We also show that a H\"older function of order $\a$, $0<\a<1$, that is analytic in the upper half-plane must be operator H\"older of order $\a$. Then we generalize these results to higher order operator differences. We obtain sharp conditions for the existence of operator derivatives and express operator derivatives in terms of multiple operator integrals with respect to semi-spectral measures. 
Finally, we obtain sharp estimates in the case of perturbations of Schatten-von Neumann class $\bS_p$ and obtain analogs of all the results for commutators and quasicommutators. Note that the proofs in the case of dissipative operators are considerably more complicated than the proofs of the corresponding results for self-adjoint operators, unitary operators, and contractions that were obtained earlier in \cite{AP2}, \cite{AP3}, and \cite{Pe6}.
\end{abstract}

\maketitle

\

\begin{center}
{\Large Contents}
\end{center}

\footnotesize

\

\begin{enumerate}
\item[1.] Introduction \quad\dotfill \pageref{intr}
\item[2.] Function spaces \quad\dotfill \pageref{bes}
\item[3.] Multiple operator integrals  \quad\dotfill \pageref{koi}
\item[4.] Dissipative operators  \quad\dotfill \pageref{dis}
\item[5.] Operator Lipschitz functions and operator differentiability \quad\dotfill \pageref{OLD}
\item[6.] Hilbert--Schmidt perturbations \quad\dotfill \pageref{HS}
\item[7.] H\"older classes and general moduli of continuity \quad\dotfill \pageref{MC}
\item[8.] Higher order operator differences \quad\dotfill \pageref{HZ}
\item[9.] Higher operator derivatives \quad\dotfill \pageref{dif}
\item[10.] Estimates in Schatten--von Neumann classes \quad\dotfill \pageref{SvN}
\item[11.] Commutators and quasicommutators \quad\dotfill \pageref{cq}
\item[] References \quad\dotfill \pageref{bibl}
\end{enumerate}

\normalsize

\

\setcounter{section}{0}
\section{\bf Introduction}
\setcounter{equation}{0}
\label{intr}

\

It is well known that a Lipschitz function on the real line is not necessarily {\it operator Lipschitz}, i.e.,
the condition 
$$
|f(x)-f(y)|\le\const|x-y|,\quad x,~y\in\R,
$$
does not imply that for self-adjoint operators
$A$ and $B$ on Hilbert space,
$$
\|f(A)-f(B)\|\le\const\|A-B\|.
$$
The existence of such functions was proved in \cite{F1}. 
Then it was shown in \cite{Ka} that the function $f(x)=|x|$ is not operator Lipschitz. Note that earlier it was shown in \cite{Mc} is not commutator Lipschitz.
Later in \cite{Pe1} necessary conditions were found for a function $f$ to be operator Lipschitz.
Those necessary conditions also imply that Lipschitz  functions do not have to be operator Lipschitz.
In particular, it was shown in \cite{Pe1} that an operator Lipschitz function must belong locally to the Besov space $B_1^1(\R)$ (see \S\,2 for an introduction to Besov spaces). Note that in \cite{Pe1} and \cite{Pe3} a stronger necessary condition was also obtained. 

It is also well known that the fact that $f$ is continuously differentiable does not imply that for bounded self-adjoint operators $A$ and $K$ the function
$$
t\mapsto f(A+tK)
$$
is differentiable. For this map to be differentiable, $f$ it must satisfy locally the same necessary conditions \cite{Pe1}, \cite{Pe3}.

On the other hand, it was proved in \cite{Pe1} and \cite{Pe3} that the condition that a function belongs to the Besov space $B_{\be1}^1(\R)$ is sufficient for operator Lipschitzness (as well as for operator differentiability). 

It turned out, however, that the situation changes  dramatically if we consider H\"older classes $\L_\a(\R)$ with $0<\a<1$. 
It was shown in \cite{AP1} and \cite{AP2} that H\"older functions are necessarily  operator H\"older, i.e., the condition
\bay
\label{hol}
|f(x)-f(y)|\le\const|x-y|^\a,\quad x,~y\in\R,
\ey
implies that for self-adjoint operators $A$ and $B$ on Hilbert space,
\bay
\label{ohol}
\|f(A)-f(B)\|\le\const(1-\a)^{-1}\|f\|_{\L_\a(\R)}\|A-B\|^\a.
\ey
Note a different proof with the factor $(1-\a)^{-2}$ instead of $(1-\a)^{-1}$ was obtained in \cite{FN}.

Analogous results were obtained in \cite{AP2} (see also \cite{AP1}) for the Zygmund class $\L_1(\R)$ and for the whole scale of H\"older--Zygmund classes $\L_\a(\R)$, $0<\a<\be$; see \S\,\ref{bes} for the definition. 

We would like to mention here that in \cite{AP1}, \cite{AP2}, and \cite{AP3} we also obtained similar estimates in the case of functions in spaces $\L_\o$ and $\L_{\o,m}$ (see \S\,\ref{bes} for the definition) as well as sharp estimates in the case when the perturbation belongs to Schatten-von Neumann classes $\bS_p$. 

Note that in \cite{AP1}, \cite{AP2}, and \cite{AP3} analogs of the above results for unitary operators and contractions were also obtained.

In this paper we deal with functions of perturbed dissipative operators and we obtain analogs of the above results for maximal dissipative operators. 
In particular, we improve results of \cite{Nab}.

We would like to mention that the case of dissipative operators is considerably more complicated. In particular, the problem is that even if we deal with bounded nonself-adjoint dissipative operators, their resolvent self-adjoint dilations are necessarily unbounded. As in the case of contractions (see \cite{Pe6}, \cite{AP1}, \cite{AP2}, and \cite{AP3}), our techniques are based on double operator integrals and multiple operator integrals with respect to semi-spectral measures. However, to obtain a representation of operator differences and operator derivatives in terms of multiple operator integrals is much more difficult than in the case of contractions. 

In \S\,\ref{OLD} we obtain sharp conditions for functions analytic in the upper half-plane to be operator Lipschitz and operator differentiable.

In \S\,\ref{HS} we consider the case of Hilbert--Schmidt perturbations and characterize so-called Hilbert--Schmidt Lipschitz functions.

Section \ref{MC} is devoted to estimates for H\"older functions and for functions of classes $\L_\o(\R)$.

We obtain in \S\,\ref{HZ} estimates for higher order operator differences for H\"older--Zygmund classes and for the spaces $\L_{\o,m}(\R)$.

We study in \S\,\ref{dif} conditions, under which higher operator derivatives exist and express higher operator derivatives in terms of multiple operator integrals.

In \S\,\ref{SvN} we obtain Schatten--von Neumann estimates in the case when the perturbation belongs to Schatten-von Neumann classes. Such results can be generalized to more general ideals of operators on Hilbert space.

Finally, in \S\,\ref{cq} we obtain sharp estimates for quasicommutators $f(L)R-Rf(M)$ for maximal dissipative operators $L$ and $M$ and a bounded operator $R$ in terms of $LR-RM$.

In \S\,\ref{bes} we collect necessary information on Besov classes (and in particular, the H\"older--Zygmund classes), and spaces $\L_\o(\R)$ and $\L_{\o,m}(\R)$. In \S\,\ref{koi} we give a brief introduction to double and multiple operator integrals. An introduction to dissipative operators is given in \S\,\ref{dis}.

\

\section{\bf Function spaces}
\setcounter{equation}{0}
\label{bes}

\

{\bf 2.1. Besov classes.}
The purpose of this subsection is to give a brief introduction to Besov spaces that play an important role in problems of perturbation theory.

Let $0<p,\,q\le\be$ and $s>0$. The homogeneous Besov class $B^s_{pq}(\R)$ of functions (or distributions) on $\R$ can be defined in the following way. Let $w$ be an infinitely differentiable function on $\R$ such
that
\bay
\label{w}
w\ge0,\quad\supp w\subset\left[\frac12,2\right],\quad\mbox{and} \quad w(x)=1-w\left(\frac x2\right)\quad\mbox{for}\quad x\in[1,2].
\ey

We define the functions $W_n$ and $W^\sharp_n$ on $\R$ by
$$
\big(\F W_n\big)(x)=w\left(\frac{x}{2^n}\right),\quad
\big(\F W^\sharp_n\big)(x)=\big(\F W_n\big)(-x),\quad n\in\Z,
$$
where $\F$ is the {\it Fourier transform}:
$$
\big(\F f\big)(t)=\int_\R f(x)e^{-{\rm i}xt}\,dx,\quad f\in L^1.
$$

With every tempered distribution $f\in{\mathscr S}^\prime(\R)$ we
associate a sequences $\{f_n\}_{n\in\Z}$,
\bay
\label{wn}
f_n\df f*W_n+f*W_n^\sharp.
\ey
Initially we define the (homogeneous) Besov class $\dot B^s_{pq}(\R)$ as the set of all $f\in{\mathscr S}^\prime(\R)$
such that
\bay
\label{Wn}
\{2^{ns}\|f_n\|_{L^p}\}_{n\in\Z}\in\ell^q(\Z).
\ey
According to this definition, the space $\dot B^s_{pq}(\R)$ contains all polynomials. Moreover, the distribution $f$ is defined by the sequence $\{f_n\}_{n\in\Z}$
uniquely up to a polynomial. It is easy to see that the series $\sum_{n\ge0}f_n$ converges in ${\mathscr S}^\prime(\R)$.
However, the series $\sum_{n<0}f_n$ can diverge in general. It is easy to prove that the
series $\sum_{n<0}f_n^{(r)}$ converges uniformly on $\R$ for each nonnegative integer
$r>s-1/p$ if $q>1$ and the series $\sum_{n<0}f_n^{(r)}$ converges uniformly, whenever $r\ge s-1/p$ if $q\le1$.

Now we can define the modified (homogeneous) Besov class $B^s_{pq}(\R)$. We say that a distribution $f$
belongs to $B^s_{pq}(\R)$ if $\{2^{ns}\|f_n\|_{L^p}\}_{n\in\Z}\in\ell^q(\Z)$ and
$f^{(r)}=\sum_{n\in\Z}f_n^{(r)}$ in the space ${\mathscr S}^\prime(\R)$, where
$r$ is the minimal nonnegative integer such that $r>s-1/p$ in the case $q>1$ and
$r$ is the minimal nonnegative integer such that $r\ge s-1/p$ in the case $q\le1$.
Now the function $f$ is determined uniquely by the sequence $\{f_n\}_{n\in\Z}$ up
to a polynomial of degree less that $r$, and a polynomial $\varphi$ belongs to $B^s_{pq}(\R)$
if and only if $\deg\varphi<r$.

It is known that the H\"older--Zygmund classes $\L_\a(\R)\df B^\a_{\be}(\R)$, $\a>0$, can be described
as the classes of continuous functions $f$ on $\R$ such that
$$
\big|(\D^m_tf)(x)\big|\le\const|t|^\a,\quad t\in\R,
$$
where the difference operator $\D_t$ is defined by
$$
(\D_tf)(x)=f(x+t)-f(x),\quad x\in\R,
$$
and $m$ is the integer such that $m-1\le\a<m$.
We can introduce the following equivalent seminorm on $\L_\a(\R)$:
$$
\sup_{n\in\Z}2^{n\a}\|f_n\|_{L^\be},\quad f\in\L_\a(\R).
$$

In this paper we deal mainly with the analytic Besov classes.
Denote by ${\mathscr S}^\prime_+(\R)$ the set of all $f\in{\mathscr S}^\prime(\R)$
such that $\supp\F f\subset[0,\be)$.
We define the analytic Besov class $\big(B_{pq}^s(\R)\big)_+$ as the intersection
$B_{pq}^s(\R)$ with ${\mathscr S}^\prime_+(\R)$.
Put $\big(\L_\a(\R)\big)_+\df\L_\a(\R)\cap{\mathscr S}^\prime_+(\R)$.
It should be noted that in the analytic case the formula \rf{wn}
can be simplified as follows
$$
f_n=f*W_n
$$
because $f*W_n^\sharp=0$ for all $f\in{\mathscr S}^\prime_+(\R)$ and $n\in\Z$.

The functions in the analytic Besov spaces admits a natural continuation to
the upper half-plane $\C_+\df\{z\in\C:\im z>0\}$. This continuation is analytic in
$\C_+$. Thus, we can consider the analytic Besov classes as spaces of function
analytic in $\C_+$. It is known that the class $\big(B_{pq}^s(\R)\big)_+$
is (in general, up to polynomials) the space of all functions $f$ analytic in $\C_+$ and such that
$$
\int_0^\be y^{\,q(m-s)-1}\left(\int_\R|f^{(m)}(x+{\rm i}y)|^p\,dx\right)^{\frac qp}\,dy<\be
$$
for some $m>s$
(with the natural modification in the case where $p=\be$ or $q=\be$).
In particular, $\big(\L_\a(\R)\big)_+$ is the space of all functions $f$
analytic in $\C_+$ and such that
$$
\sup_{y>0}y^{m-\a}\big|f^{(m)}(x+{\rm i}y)\big|<\be,
$$
where $m\in\Z$ with $m>\a$.
To optimize the set of polynomials in $\big(\L_\a(\R)\big)_+$, we can assume that $m-1\le\a$.

We refer the reader to \cite{Pee}, \cite{T}, and \cite{Pe4} for more detailed information on Besov spaces.

To define a regularized de la Vall\'ee Poussin type kernel $V_n$, we define the $C^\be$ function $v$ on $\R$ by
\bay
\label{VP}
v(x)=1\quad\mbox{for}\quad x\in[-1,1]\quad\mbox{and}\quad v(x)=w(|x|)\quad\mbox{if}\quad |x|\ge1,
\ey
where $w$ is a function described in \rf{w}.
We define the de la Vall\'ee Poussin type functions $V_n$, $n\in\Z$, (associated with $w$) by
\bay
\label{dlvp}
\F V_n(x)=v\left(\frac{x}{2^n}\right),
\ey
where $v$ is the function given by \rf{VP}.

\medskip

{\bf 2.2. Spaces $\bs{\L_\o}$ and $\bs{\L_{\o,m}}$.} Let $\o$ be a modulus of continuity, i.e., $\o$ is a nondecreasing continuous function on $[0,\be)$
such that $\o(0)=0$, $\o(x)>0$ for $x>0$, and
$$
\o(x+y)\le\o(x)+\o(y),\quad x,~y\in[0,\be).
$$
We denote by $\L_\o(\R)$ the space of functions on $\R$ such that
$$
\|f\|_{\L_\o(\R)}\df\sup_{x\ne y}\frac{|f(x)-f(y)|}{\o(|x-y|)}.
$$
Put $\big(\L_\o(\R)\big)_+\df\L_\o(\R)\cap{\mathscr S}^\prime_+(\R)$.

Consider now moduli of continuity of higher order. For a continuous function $f$ on $\R$, we define the $m$th modulus of continuity $\o_{f,m}$ of $f$ by
$$
\o_{f,m}(x)=\sup_{\{h:0\le h\le x\}}\big\|\D_h^mf\big\|_{L^\be}=\sup_{\{h:0\le|h|\le x\}}\big\|\D_h^mf\big\|_{L^\be},\quad x>0.
$$
It is well known that
$\o_{f,m}(2x)\le 2^m\o_{f,m}(x)$, see, e.g.,, \cite{DL}.

Suppose now that
$\o$ is a nondecreasing function on $(0,\be)$ such that
\bay
\label{on}
\lim_{x\to0}\o(x)=0\quad\mbox{and}\quad
\o(2x)\le2^m\o(x)\quad\mbox{for}\quad x>0.
\ey
Denote by $\L_{\o,m}(\R)$ the set of continuous functions $f$ on $\R$
satisfying
$$
\|f\|_{\L_{\o,m}(\R)}\df\sup\limits_{t>0}\frac{\|\D^m_tf\|_{L^\infty}}{\o(t)}<+\infty.
$$
Put $\big(\L_{\o,m}(\R)\big)_+\df\L_{\o,m}(\R)\cap{\mathscr S}^\prime_+(\R)$.

Note that the spaces $\L_\o(\R)$ and $\big(\L_\o(\R)\big)_+$ are special cases
of the spaces $\L_{\o,m}(\R)$ and $\big(\L_{\o,m}(\R)\big)_+$ that correspond to $m=1$.

It can be verified that a function $f$ in $\L_{\o,m}(\R)$ belongs
to the space $\big(\L_{\o,m}(\R)\big)_+$ if and only if
it has a (unique) continuous extension to the closed upper half-plane 
$\clos\C_+$ that
is analytic in the open upper half-plane $\C_+$ with at most a polynomial growth rate at infinity. We use the same notation $f$ for its extension.

We need the following inequalities
$$
\|f-f*V_n\|_{L^\be}\le\const\,\o\big(2^{-n}\big)\|f\|_{\L_{\o,m}(\R)},\quad n\in\Z,
$$
and
$$
\|f*W_n\|_{L^\be}\le\const\o\big(2^{-n}\big)\|f\|_{\L{_{\o,m}}(\R)},\quad n\in\Z,
$$
for all $f\in\L_{\o,m}(\R)$, where $\o$ is a nondecreasing function on $(0,\be)$ satisfying \rf{on},
see, for example, Theorem 2.6 and Corollary 2.7 in \cite{AP2}.

\

\section{\bf Multiple operator integrals}
\setcounter{equation}{0}
\label{koi}

\

{\bf 3.1. Double operator integrals.}
In this subsection we give a brief introduction into the theory of double operator integrals
developed by Birman and Solomyak in \cite{BS1}, \cite{BS2}, and \cite{BS3}, see also their survey \cite{BS5}.

Let $(\X,E_1)$ and $(\Y,E_2)$ be spaces with spectral measures $E_1$ and $E_2$
on Hilbert space. Let us first define double operator integrals
\bay
\label{doi}
\int\limits_{\X}\int\limits_{\Y}\Phi(x,y)\,d E_1(x)\,Q\,dE_2(y),
\ey
for bounded measurable functions $\Phi$ and operators $Q$
of Hilbert--Schmidt class $\bS_2$. Consider the set function $F$ defined 
on measurable rectangles by
$$
F(\D_1\times\D_2)Q=E_1(\D_1)QE_2(\D_2),\quad Q\in\bS_2,
$$ 
$\D_1$ and $\D_2$ being measurable subsets of $\X$ and $\Y$. Clearly, the values of $F$ are orthogonal projections on the Hilbert space $\bS_2$.

 It was shown in \cite{BS4} that $F$ extends to a spectral measure on 
$\X\times\Y$. If $\Phi$ is a bounded measurable function on $\X\times\Y$, we define
$$
\int\limits_{\X}\int\limits_{\Y}\Phi(x,y)\,d E_1(x)\,Q\,dE_2(y)=
\left(\,\,\int\limits_{\X_1\times\X_2}\Phi\,dF\right)Q.
$$
Clearly,
$$
\left\|\,\,\int\limits_{\X}\int\limits_{\Y}\Phi(x,y)\,dE_1(x)\,Q\,dE_2(y)\right\|_{\bS_2}
\le\|\Phi\|_{L^\be}\|Q\|_{\bS_2}.
$$


If the transformer
$$
Q\mapsto\int\limits_{\X}\int\limits_{\Y}\Phi(x,y)\,d E_1(x)\,Q\,dE_2(y)
$$
maps the trace class $\bS_1$ into itself, we say that $\Phi$ is a {\it Schur multiplier of $\bS_1$ associated with 
the spectral measures $E_1$ and $E_2$}. In
this case the transformer
\bay
\label{tra}
Q\mapsto\int\limits_{\Y}\int\limits_{\X}\Phi(x,y)\,d E_2(y)\,Q\,dE_1(x)
\ey
extends by duality to a bounded linear transformer on the space of bounded linear operators
and we say that the function $\Psi$ on $\X_2\times\X_1$ defined by 
$$
\Psi(y,x)=\Phi(x,y)
$$
is {\it a Schur multiplier of the space of bounded linear operators associated with $E_2$ and $E_1$}.
We denote the space of such Schur multipliers by $\fM(E_2,E_1)$.

Birman in Solomyak established in \cite{BS3} that if
 $A$ is a self-adjoint operator (not necessarily bounded),
$K$ is a bounded self-adjoint operator, and
$f$ is a continuously differentiable 
function on $\R$ such that the divided difference $\dg f$ defined by
$$
(\dg f)(x,y)\df\frac{f(x)-f(y)}{x-y},\quad x\ne y,\quad (\dg f)(x,x)\df f'(x)\quad x,\,y\in\R.
$$
belongs to $\fM(E_{A+K},E_A)$, then
\bay
\label{BSF}
f(A+K)-f(A)=\iint\limits_{\R\times\R}\big(\dg f\big)(x,y)\,dE_{A+K}(x)K\,dE_A(y)
\ey
and
$$
\|f(A+K)-f(A)\|\le\const\|\dg f\|_{\fM}\|K\|,
$$
where $\|\dg f\|_{\fM}$ is the norm of $\dg f$ in $\fM(E_{A+K},E_A)$.

In the case when $K$ belongs to the Hilbert Schmidt class $\bS_2$, formula \rf{BSF} was established in \cite{BS3} for all Lipschitz functions $f$ and it was shown that
$$
\|f(A+K)-f(A)\|_{\bS_2}\le\|f\|_{\Li}\|K\|_{\bS_2}.
$$
where 
$$
\|f\|_{\Li}\df\sup_{x\ne y}\frac{|f(x)-f(y)|}{|x-y|}.
$$
Note that  if $\f$ is not differentiable, $\dg \f$ is not defined on the diagonal of $\R\times\R$, but formula
\rf{BSF} still holds if we define $\dg \f$ to be zero on the diagonal.

It is easy to see that if a function $\Phi$ on $\X\times\Y$ belongs to the {\it projective tensor
product}
$L^\be(E_1)\hat\otimes L^\be(E_2)$ of $L^\be(E_1)$ and $L^\be(E_2)$ (i.e., $\Phi$ admits a representation
\bay
\label{ptp}
\Phi(x,y)=\sum_{n\ge0}\f_n(x)\psi_n(y),
\ey
where $\f_n\in L^\be(E_1)$, $\psi_n\in L^\be(E_2)$, and
\bay
\label{ptpn}
\sum_{n\ge0}\|\f_n\|_{L^\be}\|\psi_n\|_{L^\be}<\be),
\ey
then $\Phi\in\fM(E_1,E_2)$, i.e., $\Phi$ is a Schur multiplier of the space of bounded linear operators. For such functions $\Phi$ we have
$$
\int\limits_\X\int\limits_\Y\Phi(x,y)\,d E_1(x)Q\,dE_2(y)=
\sum_{n\ge0}\left(\,\int\limits_\X \f_n\,dE_1\right)Q\left(\,\int\limits_\Y \psi_n\,dE_2\right).
$$ 
For $\Phi$ in the projective tensor product $L^\be(E_1)\hat\otimes L^\be(E_2)$, its norm in $L^\be(E_1)\hat\otimes L^\be(E_2)$
is, by definition, the infimum of the  left-hand side of \rf{ptpn} over all representations \rf{ptp}.


More generally, $\Phi$ is a Schur multiplier  if $\Phi$ 
belongs to the {\it integral projective tensor product} $L^\be(E_1)\hat\otimes_{\rm i}L^\be(E_2)$ of $L^\be(E_1)$ and $L^\be(E_2)$, i.e., $\Phi$ admits a representation
$$
\Phi(x,y)=\int_\O \f(x,\o)\psi(y,\o)\,d\s(\o),
$$
where $(\O,\s)$ is a $\s$-finite measure space, $\f$ is a measurable function on $\X\times \O$,
$\psi$ is a measurable function on $\Y\times \O$, and
$$
\int_\O\|\f(\cdot,\o)\|_{L^\be(E_1)}\|\psi(\cdot,\o)\|_{L^\be(E_2)}\,d\s(\o)<\be.
$$
If $\Phi\in L^\be(E_1)\hat\otimes_{\rm i}L^\be(E_2)$, then
\bay
\label{iptr}
\iint\limits_{\X\times\Y}\!\Phi(x,y)\,d E_1(x)\,Q\,dE_2(y)\!=\!\!
\int\limits_\O\!\left(\,\int\limits_\X \f(x,\o)\,dE_1(x)\!\right)\!Q\!
\left(\,\int\limits_\Y \psi(y,\o)\,dE_2(y)\!\right)\!d\s(\o).
\ey

It turns out that all Schur multipliers of the space of bounded linear operators can be obtained in this way (see \cite{Pe1}).

In connection with the Birman--Solomyak formula, it is important to obtain sharp estimates of divided differences in integral projective tensor products of $L^\be$ spaces. It was shown  in \cite{Pe3} that if $f$ is a bounded function on $\R$ whose Fourier transform is supported on $[-\s,\s]$
(in other words, $f$ is an entire function of exponential type at most $\s$ that is bounded on $\R$), then $\dg f\in L^\be\hat\otimes_{\rm i}L^\be$
and
\bay
\label{Be}
\big\|\dg f\big\|_{L^\be\hat\otimes_{\rm i} L^\be}\le\const \s\|f\|_{L^\be(\R)}.
\ey
Inequality \rf{Be} led in \cite{Pe3} to the fact that functions in  $B_{\be1}^1(\R)$ are operator Lipschitz.

It was observed in \cite{Pe3} that it follows from \rf{BSF} and \rf{Be} that if $f$ is an entire function of exponential type at most $\s$ that is bounded on $\R$, and $A$ and $B$ are self-adjoint operators with bounded $A-B$, then
$$
\|f(A)-f(B)\|\le\const\s\|f\|_{L^\be}\|A-B\|.
$$
Actually, it turns out that the last inequality holds with constant equal to 1. This was established in \cite{AP4}.

\medskip

{\bf 3.2. Multiple operator integrals.} Formula \rf{iptr} suggests an approach to multiple operator integrals that is based on integral projective tensor products. This approach was given in \cite{Pe5}.

To simplify the notation, we consider here the case of triple operator integrals; the case of arbitrary multiple operator integrals can be treated in the same way.

Let $(\X,E_1)$, $(\Y,E_2)$, and $(\cZ,E_3)$
be spaces with spectral measures $E_1$, $E_2$, and $E_3$. Suppose that
the function $\Phi$ belongs to the integral projective tensor product
$L^\be(E_1)\hat\otimes_{\rm i}L^\be(E_2)\hat\otimes_{\rm i}L^\be(E_3)$, i.e., $\Phi$ admits a representation
\bay
\label{ttp}
\Phi(x,y,z)=\int_\O \f(x,\o)\psi(y,\o)\chi(z,\o)\,d\s(\o),
\ey
where $(\O,\s)$ is a $\s$-finite measure space, $\f$ is a measurable function on $\X\times \O$,
$\psi$ is a measurable function on $\Y\times \O$, $\chi$ is a measurable function on $\cZ\times \O$,
and
$$
\int_\O\|\f(\cdot,\o)\|_{L^\be(E)}\|\psi(\cdot,\o)\|_{L^\be(F)}\|\chi(\cdot,\o)\|_{L^\be(G)}\,d\s(\o)<\be.
$$

Suppose now that $T_1$  and $T_2$ are a bounded linear operators. For a function $\Phi$ in
$L^\be(E_1)\hat\otimes_{\rm i}L^\be(E_2)\hat\otimes_{\rm i}L^\be(E_3)$ of the form \rf{ttp}, we put
\begin{align}
\label{opr}
&\int\limits_\X\int\limits_\Y\int\limits_\cZ\Phi(x,y,z)
\,d E_1(x)T_1\,dE_2(y)T_2\,dE_3(z)\\[.2cm]
\df&\int\limits_\O\left(\,\int\limits_\X \f(x,\o)\,dE_1(x)\right)T_1
\left(\,\int\limits_\Y \psi(y,\o)\,dE_2(y)\right)T_2
\left(\,\int\limits_\cZ \chi(z,\o)\,dE_3(z)\right)\,d\s(\o).\nonumber
\end{align}

It was shown in \cite{Pe5} (see also \cite{ACDS} for a different proof)  that the above definition does not depend on the choice of a representation \rf{ttp}.
%

It is easy to see that the following inequality holds
$$
\left\|\int\limits_\X\int\limits_\Y\int\limits_\cZ\Phi(x,y,z)
\,dE_1(x)T_1\,dE_2(y)T_2\,dE_3(z)\right\|
\le\|\Phi\|_{L^\be\hat\otimes_{\rm i}L^\be\hat\otimes_{\rm i}L^\be}\cdot\|T_1\|\cdot\|T_2\|.
$$

In particular, the triple operator integral on the left-hand side of \rf{opr} can be defined if $\Phi$ belongs to the projective
tensor product $L^\be(E_1)\hat\otimes L^\be(E_2)\hat\otimes L^\be(E_3)$, i.e., $\Phi$ admits a representation
$$
\Phi(x,y,z)=\sum_{n\ge1}\f_n(x)\psi_n(y)\chi_n(z),
$$
where $\f_n\in L^\be(E_1)$, $\psi_n\in L^\be(E_2)$, $\chi_n\in L^\be(E_3)$ and
$$
\sum_{n\ge1}\|\f_n\|_{L^\be(E_1)}\|\psi_n\|_{L^\be(E_2)}\|\chi_n\|_{L^\be(E_3)}<\be.
$$

In a similar way one can define multiple operator integrals, see \cite{Pe5}.

Recall that multiple operator integrals were considered earlier in \cite{Pa} and \cite{S}. However, in those papers the class of functions 
$\Phi$ for which the left-hand side of \rf{opr} was defined is much narrower than in the definition given above.

Multiple operator integrals are used in \cite{Pe5} in connection with the problem of evaluating higher order operator derivatives. 
To obtain formulae for higher operator derivatives, one has to integrate divided differences of higher orders (see \cite{Pe5}). 

For a function $f$ on $\R$, the {\it divided differences $\dg^m f$ of order $m$} are defined inductively as follows:
$$
\dg^0f\df f;
$$
if $m\ge1$, then in the case when $x_1,x_2,\cdots,x_{m+1}$ are distinct points in $\R$,
$$
(\dg^{m}f)(x_1,\cdots,x_{m+1})\df
\frac{(\dg^{m-1}f)(x_1,\cdots,x_{m-1},x_m)-
(\dg^{m-1}f)(x_1,\cdots,x_{m-1},x_{m+1})}{x_{m}-x_{m+1}}
$$
(the definition does not depend on the order of the variables). Clearly,
$$
\dg f=\dg^1f.
$$
If $f\in C^m(\R)$, then $\dg^{m}f$ extends by continuity to a function defined for all points $x_1,x_2,\cdots,x_{m+1}$.

It can be shown that
$$
({\frak D}^m f)(x_1,\dots,x_{m+1})=\sum\limits_{k=1}^{m+1}f(x_k)
\prod\limits_{j=1}^{k-1}(x_k-x_j)^{-1}\prod\limits_{j=k+1}^{m+1}(x_k-x_j)^{-1}.
$$
 
It was established in \cite{Pe5} that if $f$ is an entire function of exponential type at most $\s$ that is bounded on $\R$, then
\bay
\label{Boke}
\big\|\dg^m f\big\|_{L^\be\hat\otimes_{\rm i}\cdots\hat\otimes_{\rm i} L^\be}\le\const \s^m\|f\|_{L^\be(\R)}.
\ey

Note that recently in \cite{JTT} Haagerup tensor products were used to define multiple operator integrals. 

\medskip

{\bf 3.3. Multiple operator integrals with respect to semi-spectral measures.}
Let $\h$ be a Hilbert space and let $(\X,{\frak B})$ be a measurable space.
A map $\E$ from ${\frak B}$ to the algebra ${\mathscr B}(\h)$ of all bounded operators on $\h$ is called a {\it semi-spectral measure}
if 
$$
\E(\D)\ge\0,\quad\D\in{\frak B},
$$
$$
\E(\varnothing)=\0\quad\mbox{and}\quad\E(\X)=I,
$$
and for a sequence $\{\D_j\}_{j\ge1}$ of disjoint sets in ${\frak B}$,
$$
\E\left(\bigcup_{j=1}^\be\D_j\right)=\lim_{N\to\be}\sum_{j=1}^N\E(\D_j)\quad\mbox{in the weak operator topology}.
$$

\medskip

If $\K$ is a Hilbert space such that $\h\subset\K$ and $E:{\frak B}\to{\mathscr B}(\K)$ is a spectral measure on $(\X,{\frak B})$, then it is easy to see that the map $\E:{\frak B}\to {\mathscr B}(\h)$ defined by
\bay
\label{dil}
\E(\D)=P_\h E(\D)\big|\h,\quad\D\in{\frak B},
\ey
is a semi-spectral measure. Here $P_\h$ stands for the orthogonal projection onto $\h$.

Naimark proved in \cite{N}  that all semi-spectral measures can be obtained in this way, i.e.,
a semi-spectral measure is always a {\it compression} of a spectral measure. A spectral measure $E$ satisfying \rf{dil} is called a {\it spectral dilation of the semi-spectral measure} $\E$.

A spectral dilation $E$ of a semi-spectral measure $\E$ is called {\it minimal} if 
$$
\K=\clos\spn\{E(\D)\h:~\D\in{\frak B}\}.
$$

It was shown in \cite{MM} that if $E$ is a minimal spectral dilation of a semi-spectral measure $\E$, then
$E$ and $\E$ are mutually absolutely continuous and all minimal spectral dilations of a semi-spectral measure are isomorphic in the natural sense.

If $\f$ is a bounded complex-valued measurable function on $\X$ and $\E:{\frak B}\to {\mathscr B}(\h)$ is a semi-spectral measure, then the integral
\bay
\label{iss}
\int_\X \f(x)\,d\E(x)
\ey
can be defined as
\bay
\label{voi}
\int_\X \f(x)\,d\E(x)=\left.P_\h\left(\int_\X \f(x)\,d E(x)\right)\right|\h,
\ey
where $E$ is a spectral dilation of $\E$. It is easy to see that the right-hand side of \rf{voi} does not depend on the choice
of a spectral dilation. The integral \rf{iss} can also be computed as the limit of sums
$$
\sum \f(x_\a)\E(\D_\a),\quad x_\a\in\D_\a,
$$
over all finite measurable partitions $\{\D_\a\}_\a$ of $\X$.

If $T$ is a contraction on a Hilbert space $\h$, then by the Sz.-Nagy dilation theorem
(see \cite{SNF}),  $T$ has a unitary dilation, i.e., there exist a Hilbert space $\K$ such that
$\h\subset\K$ and a unitary operator $U$ on $\K$ such that
\bay
\label{DT}
T^n=P_\h U^n\big|\h,\quad n\ge0,
\ey
where $P_\h$ is the orthogonal projection onto $\h$. Let $E_U$ be the spectral measure of $U$.
Consider the operator set function $\E$ defined on the Borel subsets of the unit circle $\T$ by
$$
\E(\D)=P_\h E_U(\D)\big|\h,\quad\D\subset\T.
$$
Then $\E$ is a semi-spectral measure. It follows immediately from
\rf{DT} that 
\bay
\label{step}
T^n=\int_\T \z^n\,d\E(\z)=P_\h\int_\T\z^n\,dE_U(\z)\Big|\h,\quad n\ge0.
\ey
It is easy to see that $\E$ does not depend on the choice of a unitary dilation.
$\E$ is called the {\it semi-spectral measure} of $T$.

It follows easily from  \rf{step} that 
$$
f(T)=P_\h\int_\T f(\z)\,dE_U(\z)\Big|\h
$$
for an arbitrary function $f$ in the disk-algebra $C_A$.

In \cite{Pe2} and  \cite{Pe6} double operator integrals and multiple operator integrals with respect to semi-spectral measures were introduced.

For semi-spectral measures $\E_1$ and $\E_2$, double operator integrals
$$
\iint\limits_{\X_1\times\X_2}\Phi(x_1,x_2)\,d\E_1(x_1)Q\,d\E_2(X_2).
$$
were defined in \cite{Pe6} by analogy with the case of spectral measures in the case when $Q\in\bS_2$ and $\Phi$ is a bounded measurable function and in the case when $Q$ is a bounded linear operator and $\Phi$ belongs to the integral projective tensor product of the spaces $L^\be(\E_1)$ and $L^\be(\E_2)$.

Similarly, multiple operator integrals with respect to semi-spectral measures were defined in \cite{Pe6} for functions that belong to the integral projective tensor product of the corresponding $L^\be$ spaces.

\

\section{\bf Dissipative operators}
\setcounter{equation}{0}
\label{dis}

\

This section is a brief introduction into the theory of dissipative operators.
We refer the reader to \cite{SNF} and \cite{So} for more information.

\medskip

{\bf Definition.}
Let $\h$ be a Hilbert space. An operator $L$ (not necessarily bounded)
with dense domain $\cd_L$ in $\h$ is called {\it dissipative} if
$$
\im(Lu,u)\ge0,\quad u\in\cd_L.
$$
A dissipative operator is called {\it maximal dissipative} if it has no proper dissipative extension.

\medskip

Note that if $L$ is a symmetric operator (i.e., $(Lu,u)\in\R$ for every
$u\in\cd_L$), then $L$ is dissipative. However, it can happen that $L$ is maximal symmetric, but not maximal dissipative.

The {\it Cayley transform} of a dissipative operator $L$ is defined by
$$
T\df(L-{\rm i}I)(L+{\rm i}I)^{-1}
$$
with domain $\cd_T=(L+{\rm i}I)\cd_L$ and range $\Range T=(L-{\rm i}I)\cd_L$
(the operator $T$ is not densely defined in general). $T$ is a contraction, i.e., $\|Tu\|\le\|u\|$, $u\in\cd_T$, $1$ is not an eigenvalue of $T$, and $\Range(I-T)\df\{u-Tu:~u\in\cd_T\}$ is dense.

Conversely, if $T$ is a contraction defined on its domain $\cd_T$, $1$ is not an eigenvalue of $T$, and $\Range(I-T)$ is dense, then it is the Cayley transform of a dissipative operator $L$ and $L$ is the inverse Cayley transform of $T$:
$$
L={\rm i}(I+T)(I-T)^{-1},\quad \cd_L=\Range(I-T).
$$

A dissipative operator is maximal if and only if the domain of its Cayley transform is the whole Hilbert space.

Every dissipative operator $L$ has a maximal dissipative extension.
Every maximal dissipative operator is necessarily closed.

If $L$ is a maximal dissipative operator, then $-L^*$ is also maximal dissipative. 

If $L$ is a maximal dissipative operator, then its spectrum $\s(L)$ is contained in the upper half-plane $\C_+$ and
\bay
\label{res}
\big\|(L-\l I)^{-1}\big\|\le\frac1{|\im\l|},\quad\im\l<0.
\ey

If $L$ and $M$ are maximal dissipative operators, we say that the {\it operator $L-M$ is bounded} if there exists a bounded operator $K$ such that $M=L+K$.

We will need the following elementary fact:

\begin{lem}
\label{ef}
Let $L$ be a maximal dissipative operator and let $M$ be a dissipative operator
such that $\cd(L)=\cd(M)$ and $L-M$ is a bounded operator. Then $M$ is a maximal dissipative operator.
\end{lem}

\Pf Put $Q=L-M$ and let $\vk>\|Q\|$. Since $L$ is maximal dissipative, we have
by \rf{res},
\bay
\label{vk}
\big\|(L+{\rm i}\vk I)u\big\|\ge\vk\|u\|,\quad u\in\h,
\ey
where $\h$ is our Hilbert space. Moreover, the operator $\vk^{-1}L$ is
also maximal dissipative and
since the domain of the Cayley transform of $\vk^{-1}L$ is
$\Range(L+{\rm i}\vk I)$, it follows that
$\Range(L+{\rm i}\vk I)=\h$.

It follows easily from \rf{vk} and from the fact that $\|Q\|<\vk$ that
$\Range(L-Q+{\rm i}\vk I)=\Range(M+{\rm i}\vk I)=\h$. Thus
$\vk^{-1}M$ is maximal dissipative, and so is $M$. $\bl$

\medskip

{\bf Definition.}
A maximal dissipative operator $L$ in a Hilbert space $\h$ is called {\it pure} if there are no nonzero
(closed) subspace $\K$ of $\h$ such that $L$ induces a self-adjoint operator in $\K$.

\medskip

As it is well known (see \cite{SNF}, Prop.
4.3 of Ch. IV) for every dissipative operator $L$ in $\h$, there exists 
a unique decomposition $\h=\h_{0}\oplus\h_{\rm 1}$ reducing
$L$ and such that $L\big|\h_{0}$ is a self-adjoint
operator and $L\big|\h_{1}$ is a pure maximal dissipative operator.
(We say that a decomposition $\h=\h_{0}\oplus\h_{\rm 1}$ reduces
an operator $L$ in $\h$ if there exist (unique) operators $L_0$ in $\h_{0}$
and $L_1$ in $\h_{1}$ such that $L=L_0\oplus L_1$.)

We use the notation
\bay
\label{dec}
\h_{\rm sa}\df\h_0\quad\mbox{and}\quad\h_{\rm p}\df\h_1.
\ey

We proceed now to the construction of functional calculus for dissipative operators. Let $L$ be a maximal dissipative operator and let $T$ be its Cayley transform. Consider its minimal unitary dilation $U$, i.e., $U$ is a unitary operator defined on a Hilbert space $\K$ that contains $\h$ such that
$$
T^n=P_\h U^n\big|\h,\quad n\ge0,
$$
and $\K=\clos\spn\{U^nh:~h\in\h\}$. Since $1$ is not an eigenvalue of $T$, it follows that $1$ is not an eigenvalue of $U$ (see \cite{SNF}, Ch. II, \S\,6).

The Sz.-Nagy--Foia\c s functional calculus allows us to define a functional calculus for $T$ on the Banach algebra
$$
C_{A,1}\df\big\{g\in H^\be:~g\quad
\mbox{is continuous on}\quad\T\setminus\{1\}~\big\}.
$$
If $g\in C_{A,1}$, we put
$$
g(T)\df P_\h g(U)\Big|\h.
$$
This functional calculus is linear and multiplicative and
$$
\|g(T)\|\le\|g\|_{H^\be},\quad g\in C_{A,1}.
$$

We can define now a functional calculus for our dissipative operator on the
Banach algebra
$$
C_{A,\be}=\big\{f\in H^\be(\C_+):~f\quad\mbox{is continuous on}\quad
\R\big\}.
$$
Indeed, if $f\in C_{A,\be}$, we put
$$
f(L)\df \big(f\circ\o\big)(T),
$$
where $\o$ is the conformal map of $\dd$ onto $\C_+$ defined by
$\o(\z)\df{\rm i}(1+\z)(1-\z)^{-1}$, $\z\in\dd$.

The reader can find more detailed information in \cite{SNF}.

We proceed now to the definition of a resolvent self-adjoint dilation of a maximal dissipative operator. If $L$ is a maximal dissipative operator on a Hilbert space $\h$, we say that a self-adjoint operator $A$ in a Hilbert space
$\K$, $\K\supset\h$, is called a {\it resolvent self-adjoint dilation} of $L$
if
$$
(L-\l I)^{-1}=P_\h(A-\l I)^{-1}\Big|\h,\quad \im\l<0.
$$
The dilation is called {\it minimal} if
$$
\K=\clos\spn\big\{(A-\l I)^{-1}v:~v\in\h,~\im\l<0\big\}.
$$
If $f\in C_{A,\be}$, then
$$
f(L)=P_\h f(A)\Big|\h,\quad f\in C_{A,\be}.
$$

A minimal resolvent self-adjoint dilation of a maximal dissipative operator always exists (and is unique up to a natural isomorphism). Indeed, it suffices to take a minimal unitary dilation of the Cayley transform of this operator and
apply the inverse Cayley transform to the minimal unitary dilation.

Sometimes mathematicians use the term "self-adjoint dilation" rather than
"resolvent self-adjoint dilation". However, we believe that the term "self-adjoint dilation" is misleading.

Let us define now the semi-spectral measure of a maximal dissipative operator $L$. Let $T$ be the Cayley transform of $L$ and let $\E_T$ be the semi-spectral measure of $T$ on the unit circle $\T$. Then
\bay
\label{ssu}
g(T)=\int_\T g(\z)\,d\E_T(\z),\quad g\in C_{A,1}.
\ey
We can define now the semi-spectral measure $\E_L$ of $L$ by
$$
\E_L(\D)=\E_T\big(\o^{-1}(\D)\big),\quad\D\quad\mbox{is a Borel subset of}
\quad\R.
$$
It follows easily from \rf{ssu} that
$$
f(L)=\int_\R f(x)\,d\E_L(x),\quad f\in C_{A,\be}.
$$

We conclude this section with the following well-known lemma.

\begin{lem}
\label{61}
Let $\E_L$ be the semi-spectral measure of a maximal dissipative operator
$L$ in a Hilbert space $\h$.
Then  $\h_{\rm p}$ and
$\h_{\rm sa}$ (defined by {\em\rf{dec}}) are invariant subspaces of $\E_L$, the restriction of $\E_L$
to $\h_{\rm p}$ is an absolutely continuous (with respect to Lebesque
measure) semi-spectral measure, while the restriction of $\E_L$ to
$\h_{\rm sa}$ is a spectral measure.
\end{lem}

\Pf It suffices to prove the corresponding result for the Cayley transform $T$ of $L$
and use the well-known fact that the minimal unitary dilation of the
completely nonunitary contraction $T\big|\h_{\rm p}$
is a unitary operator with absolutely continuous spectrum, see Th. 6.4
of Ch. II
in \cite{SNF}. $\bl$

\

\section{\bf Operator Lipschitz functions and operator differentiability}
\setcounter{equation}{0}
\label{OLD}

\

In this section we estimate the norm of $f(L)-f(M)$ for an entire function $f$ of exponential type $\s$ bounded on $\R$ and for maximal dissipative operators $L$ and $M$ in terms of $\|L-M\|$. We express $f(L)-f(M)$ in terms of double operator integrals. This allows us to prove that functions in $B_{\be1}^1(\R)$ are operator Lipschitz and operator differentiable and we express the operator derivative
$$
\frac{d}{dt}f\big(L+t(M-L)\big)
$$
in terms of double operator integrals.

\begin{thm}
\label{OLB}
Let $f\in \big(B_{\be1}^1(\R)\big)_+$ and let $L$ and $M$ be maximal dissipative operators
such that $L-M$ is bounded. Then
\bay
\label{BSd}
f(L)-f(M)=\iint\frac{f(x)-f(y)}{x-y}\,d\E_L(x)(L-M)\,d\E_M(y).
\ey
\end{thm}

\begin{lem}
\label{OLe}
Let $f$ be a bounded function on $\R$ whose Fourier transform has compact support
in $(0,\be)$,
 and let $L$ and $M$ be maximal dissipative operators such that $L-M$ is bounded.
Then {\em\rf{BSd}} holds.
\end{lem}

To prove Lemma \ref{OLe}, we need  more lemmata.

\begin{lem}
\label{sot}
Let $a>0$, and
let $L$ and $M$ be maximal dissipative operators
that satisfy the hypotheses of Theorem {\em\ref{OLB}}. Then the function
$$
t\mapsto\exp({\rm i}tL)\exp\big(\ri(a-t)M\big)
$$
is differentiable in the strong operator topology on the interval $[0,a]$ and
\bay
\label{pro}
\frac{d}{dt}\exp({\rm i}tL)\exp\big(\ri(a-t)M\big)=
\ri\exp({\rm i}tL)(L-M)\exp\big(\ri(a-t)M\big).
\ey
\end{lem}

\Pf Let us observe that for an arbitrary maximal dissipative operator $L$, the
function $t\mapsto\exp(\ri tL)$, $t\ge0$, is continuous in the strong operator topology.
This follows easily from the formula
$$
\exp(\ri tL)=\int_\R e^{\ri ts}\,d\E_L(s),\quad t\ge0,
$$
where $\E_L$ is the semi-spectral measure of $L$. It follows that both functions
in \rf{pro} are continuous in the strong operator topology on $[0,a]$.

Another observation is that if $\Theta$ is an operator-valued function on an interval that is differentiable in the weak operator topology and both $\Theta$ and its derivative $\Theta'$ are continuous in the strong operator topology, then $\Theta$ is differentiable in the strong operator topology. Indeed, it suffices to represent $\Theta$ as an indefinite integral of $\Theta'$.

Thus it suffices to prove that \rf{pro} holds in the weak operator topology.

Let $u\in\cd_M=\cd_L$ and $v\in\cd_{L^*}=\cd_{M^*}$. We have
\begin{align*}
\frac{d}{dt}\Big(\exp({\rm i}tL)\exp\big(\ri(a-t)M\big)u,v\Big)&=
\frac{d}{dt}\Big(\exp\big(\ri(a-t)M\big)u,\exp(-{\rm i}tL^*)v\Big)\\[.2cm]
&=-\ri\Big(\exp\big(\ri(a-t)M\big)Mu,\exp(-{\rm i}tL^*)v\Big)\\[.2cm]
&+\ri\Big(\exp\big(\ri(a-t)M\big)u,\exp(-{\rm i}tL^*)L^*v\Big)\\[.2cm]
&=\ri\Big(\exp({\rm i}tL)(L-M)\exp\big(\ri(a-t)M\big)u,v\Big).\quad\bl
\end{align*}

\begin{cor}
\label{in}
Let $a>0$. Then
\begin{align*}
\exp(\ri aL)-\exp(\ri aM)&=\ri\int_0^a\exp({\rm i}tL)(L-M)\exp\big(\ri(a-t)M\big)\,dt\\[.2cm]
&=\ri\int_0^a\exp({\rm i}tM)(L-M)\exp\big(\ri(a-t)L\big)\,dt.
\end{align*}
\end{cor}

\begin{cor}
\label{ex}
Let $a>0$. Then
$$
\big\|\exp(\ri aL)-\exp(\ri aM)\big\|\le a\|L-M\|.
$$
\end{cor}

\begin{lem}
\label{H1}
Suppose that $L$ and $M$ satisfy the hypotheses of Theorem {\em\ref{OLB}}. Let $f$ be a function in $H^1(\C_+)$ such that
$$
\int_0^\be(1+|\xi|)\big|(\F f)(\xi)\big|\,d\xi<\be.
$$
Then
$$
f(L)-f(M)=
\frac\ri{2\pi}\int_0^\be\int_0^\be(\F f)(\xi+\eta)\exp(\ri \xi L)(L-M)\exp(\ri \eta M)\,d\xi\,d\eta.
$$
\end{lem}

\Pf Let us first observe that if $f$ satisfies the hypotheses of the Lemma, then for an arbitrary maximal dissipative operator $L$,
$$
f(L)=\frac1{2\pi}\int_0^\be(\F f)(x)\exp(\ri xL)\,dx.
$$
Indeed, it is easy to reduce this formula to the corresponding formula for the function
of a self-adjoint dilation of $L$.

We have
\begin{align*}
f(L)-f(M)&=\frac1{2\pi}\int_0^\be(\F f)(x)\big(\exp(\ri xL)-\exp(\ri xM)\big)\,dx\\[.2cm]
&=\frac\ri{2\pi}\int_0^\be(\F f)(x)\int_0^x\exp({\rm i}\xi L)(L-M)\exp\big(\ri(x-\xi)M\big)\,d\xi\,dx\\
&=\frac\ri{2\pi}\int_0^\be\int_0^\be(\F f)(\xi+\eta)\exp(\ri\xi L)(L-M)\exp(\ri\eta M)\,d\xi\,d\eta.
\quad\bl
\end{align*}


We need some results from \cite{Pe5}. For a function $f$ in $H^1(\C_+)$ and $a>0$, we define the function $f_{(a)}$ by
$$
\big(\F f_{(a)}\big)(\xi)=\frac{\xi}{\xi+a}(\F f)(\xi+a)\chi_{[0,\be)}(\xi).
$$
Put
$$
\f_{a}\df\F^{-1}\left(\min\left\{1,\frac{a}{|\xi|}\right\}\right),\quad\xi\in\R.
$$
\begin{lem} 
\label{Polya}
Let $0<a<\be$. Then $\f_a\in L^1(\R)$
and $\|\f_a\|_{L^1}\le3$.
\end{lem}
\Pf It is easy to see that $\|\f_{a}\|_{L^1}=\|\f_1\|_{L^1}$ for every $a>0$. Let us show that $\|\f_1\|_{L^1}\le3$.
We have
$$
\min\left\{1,\frac{1}{|\xi|}\right\}=\psi_1(\xi)-\psi_2(\xi),
$$
where
$$
\psi_1(\xi)\df\left((1-|\xi|)_++\min\left\{1,\frac{1}{|\xi|}\right\}\right)\qm
\psi_2(\xi)\df(1-|\xi|)_+.
$$
Clearly, both $\psi_1$ and $\psi_2$ are continuous, even and convex on $(0,\be)$. Moreover,
$\psi_1$ and $\psi_2$ vanish at $\be$.
By P\'olya's theorem (see \cite{Po}), $\F^{-1}\psi_j\in L^1$ and $\|\F^{-1}\psi_j\|_{L^1}=\psi_j(0)$, $j=1,\,2$.
Hence,
$\|\f_a\|_{L^1}\le\psi_1(0)+\psi_2(0)=3$. $\bl$

Note that
\bay
\label{defa}
f_{(a)}=e_{-a}(f-\f_a*f)=
e_{-a}f-e_{-a}(\f_{a}*f)
\ey
for every $f\in H^1(\C_+)$,
where
$$
e_\a(x)\df e^{\ri\a x},\quad x\in\R,\quad\a\in\R.
$$
The formula \rf{defa} and Lemma \ref{Polya} allow us to define $f_{(a)}$ in the case where $f\in H^\be(\C_+)$.
The following result was obtained in \cite{Pe5}. We give a proof here for completeness.
\begin{lem}
\label{fa}
Let $f\in H^\be(\C_+)$. Then $f_{(a)}\in H^\be(\C_+)$ for every $a>0$ and
\bay
\label{f_a}
\|f_{(a)}\|_{H^\be}\le4\|f\|_{H^\be}.
\ey
\end{lem}
\Pf Note that $\|f_{(a)}\|_{L^\be}\le4\|f\|_{L^\be}$ for every $f\in L^\be(\R)$,
where $f_{(a)}$ is defined by \rf{defa}.
It remains to verify that $f_{(a)}\in H^\be(\C_+)$ for every $f\in H^\be(\C_+)$.
This is clear for $f\in H^1(\C_+)\cap H^\be(\C_+)$. To complete the proof
we observe that $H^1(\C_+)\cap H^\be(\C_+)$ is dense in $H^\be(\C_+)$ in the weak-$*$
topology $\s(L^\be,L^1)$ and formula \rf{defa} defines the mapping $f\mapsto f_{(a)}$ which acts continuously
on $L^\be(\R)$ in the weak-$*$ topology $\s(L^\be,L^1)$. $\bl$

Note that it follows easily from the definition of $f_{(a)}$ that $f_{(a)}=\0$ if 
$f\in H^\be(\C_+)$ and $\supp\F f\subset(0,a)$. Hence, the same is true for all 
$f\in H^\be(\C_+)$ with $\supp\F f\subset[0,a]$.

\medskip

{\bf Proof of Lemma \ref{OLe}.} Let $\supp\F f\subset(0,\s)$, where $\s\in(0,\be)$. It was proved in \cite{Pe5} that
\bay
\label{rr}
(\dg f)(x,y)=\ri\int_0^\be e^{\ri\xi x}f_{(\xi)}(y)\,d\xi+\ri\int_0^\be f_{(\eta)}(x)e^{\ri\eta y}\,d\eta.
\ey
Since $f_{(a)}=\0$ for $a>\s$, it follows from Lemma \ref{fa} that
\rf{rr} is a representation of $\dg f$ in the integral projective tensor product
$L^\be\hat\otimes_{\rm i}L^\be$.

Thus
\begin{align*}
\iint(\dg f)(x,y)\,d\E_L(x)(L-M)\,d\E_M(y)&=\ri\int_0^\s\exp(\ri\xi L)(L-M)f_{(\xi)}(M)\,d\xi\\[.2cm]
&+\ri\int_0^\s f_{(\eta)}(L)(L-M)\exp(\ri\eta M)\,d\eta.
\end{align*}

On the other hand, if $f\in H^1(\C_+)$ and $\supp\F f\subset(0,\s)$, then
by Lemma \ref{H1},
\begin{align}
\label{mef}
f(L)&-f(M)=
\frac{\ri}{2\pi}\int_0^\be\int_0^\be(\F f)(\xi+\eta)\exp(\ri\xi L)(L-M)\exp(\ri\eta M)\,d\xi\,d\eta\nonumber
\\[.2cm]
=&\frac{\ri}{2\pi}\int_0^\be\exp(\ri\xi L)(L-M)
\left(\int_0^\be\frac{\eta(\F f)(\xi+\eta)}{\xi+\eta}\exp(\ri\eta M)\,d\eta\right)\,d\xi\nonumber\\[.2cm]
&+\frac{\ri}{2\pi}\int_0^\be\left(\int_0^\be\frac{\xi(\F f)(\xi+\eta)}{\xi+\eta}\exp(\ri\xi L)\,d\xi\right)
(L-M)\exp(\ri\eta M)\,d\eta\nonumber\\[.2cm]
=&\ri\int_0^\s\exp(\ri\xi L)(L-M)f_{(\xi)}(M)\,d\xi+
\ri\int_0^\s f_{(\eta)}(L)(L-M)\exp(\ri\eta M)\,d\eta.
\end{align}
This proves the result under the extra assumption $f\in H^1(\C_+)$.

In the general case, we approximate $f$ with the functions $f^{[\e]}$,
$$
f^{[\e]}(x)\df\dfrac{\sin^2{\e x}}{\e^2x^2}f(x).
$$
Clearly,
the support of $\mathscr F f^{[\e]}$ is contained in the  $2\e$-neighborhood
of the support $\mathscr F f$. Hence, $f^{[\e]}\in H^1(\C_+)$ and $\supp\F f^{[\e]}\subset(0,\s)$ for sufficient small $\e>0$.

Now we can apply \rf{mef} for $f^{[\e]}$ in place of $f$ and pass to the limit as $\e\to0$
in the strong operator topology.  $\bl$

\begin{lem}
\label{sle}
Let $f$ be a bounded function on $\R$ whose Fourier transform
is supported on $[0,\s]$ and let $L$ and $M$ be maximal dissipative operators
such that $L-M$ is bounded. Then
$$
\|f(L)-f(M)\|\le8\s\|f\|_{L^\be(\R)}\|L-M\|.
$$
\end{lem}

\Pf This follows from \rf{mef} and \rf{f_a}. $\bl$

\medskip

{\bf Proof of Theorem \ref{OLB}.} Let $f_n=f*W_n$, $n\in\Z$. Then
$$
\dg f=\sum_{n=-\be}^\be \dg f_n
$$
and the series converges uniformly. By Lemma \ref{OLe},
$$
f_n(L)-f_n(M)=\iint(\dg f_n)(x,y)\,d\E_L(x)(L-M)\,d\E_M(y).
$$
The result follows now from Lemma \ref{sle}. $\bl$

\begin{lem}
\label{OL}
Let $f\in\big(B_{\be1}^1\big)_+\cap H^\be(\C_+)$. Then $f$ is operator Lipschitz on the class of maximal dissipative operators.
\end{lem}

\Pf Suppose that $L$ and $M$ are maximal dissipative operators such that $L-M$ is bounded.
Let $f_n=f*W_n$, $n\in\Z$. By Lemma \ref{sle},
$$
\|f_n(L)-f_n(M)\|\le\const2^n\|f_n\|_{L^\be(\R)}\|L-M\|.
$$
It follows that
$$
\|f(L)-f(M)\|\le\const\left(\sum_{n\in\Z}2^n\|f_n\|_{L^\be}\right)\|L-M\|
\le\const\|f\|_{B^1_{\be1}(\R)}\|L-M\|.\quad\bl
$$

Note that functions in $\big(B_{\be1}^1\big)_+$ do not have to be bounded and we have not defined
unbounded functions of maximal dissipative operators. However, for functions in $\big(B_{\be1}^1\big)_+$ and for maximal dissipative operators $L$ and $M$ with bounded $L-M$, we can define the operator
$f(L)-f(M)$ by
\bay
\label{neo}
f(L)-f(M)\df\sum_{n\in\Z}\big(f_n(L)-f_n(M)\big).
\ey
It follows from Lemma \ref{sle} that the series on the right-hand side of \rf{neo} converges absolutely in the operator norm. It is also easy to verify that the right-hand side does not depend on the choice of $w$ in \rf{w}. This allows us to get rid of the condition $f\in H^\be(\C_+)$ in the statement of Theorem \ref{OL}.

\begin{thm}
\label{OLn}
Let $f\in\big(B_{\be1}^1\big)_+$. Then $f$ is operator Lipschitz on the class of maximal dissipative operators.
\end{thm}

We proceed now to operator differentiability. Suppose that $L$ and $M$ are maximal dissipative operators such that $L-M$ is bounded. Consider the family of operators
\bay
\label{Lt}
L_t=L+t(M-L),\quad t\in[0,1].
\ey
It follows from Lemma \ref{ef} that $L_t$ is maximal dissipative for every $t\in[0,1]$.
For a function $f$ of class $C_{A,\be}$ (see \S\,\ref{dis}), we consider the map
\bay
\label{fLt}
t\mapsto f(L_t)
\ey
and study conditions of its differentiability.

We are going to show that the function \rf{fLt} is differentiable in the operator norm for
$f\in\big(B_{\be1}^1\big)_+$. We have a problem similar to the problem to prove that functions in $\big(B_{\be1}^1\big)_+$ are operator Lipschitz: we have not defined unbounded functions of maximal dissipative operators. However, to differentiate the function \rf{fLt}, we consider the limit of
the operators
$$
\frac1s\big(f(L_{t+s})-f(L_t)\big),
$$
which are well defined, see \rf{neo}.

\begin{thm}
\label{der}
Let $L$, $M$ and $L_t$ be as above and let $f\in\big(B_{\be1}^1\big)_+$. Then the function
{\rm\rf{fLt}} is differentiable in the norm and
\bay
\label{DKf}
\frac{d}{ds}\big(f(L_s)-f(L)\big)\Big|_{s=t}=\iint\frac{f(x)-f(y)}{x-y}\,d\E_t(x)(M-L)\,d\E_t(y),
\ey
where $\E_t$ is the semi-spectral measure of $L_t$.
\end{thm}







{\bf Proof of Theorem \ref{der}.} As in the proof of Theorem \ref{OLB}, we put $f_n=f*W_n$. It follows from
Lemma \ref{sle} that it suffices to prove \rf{DKf} for each function $f_n$ in place of $f$. Let $g=f_n$.
For simplicity, we assume that $t=0$.

By Theorem \ref{OLB}, we have to show that
$$
\lim_{s\to0}\iint\frac{g(x)-g(y)}{x-y}\,d\E_s(x)(M-L)\,d\E_0(y)=
\iint\frac{g(x)-g(y)}{x-y}\,d\E_0(x)(M-L)\,d\E_0(y)
$$
in the norm. By \rf{rr},
\begin{align*}
\iint\frac{g(x)-g(y)}{x-y}\,d\E_s(x)(M-L)\,d\E_0(y)={\rm i}\int_0^\be e^{{\rm i}\xi L_s}g_{(\xi)}(L)\,d\xi
+{\rm i}\int_0^\be g_{(\eta)}(L_s)e^{{\rm i}\eta L}\,d\eta,
\end{align*}
while
$$
\iint\frac{g(x)-g(y)}{x-y}\,d\E_0(x)(M-L)\,d\E_0(y)={\rm i}\int_0^\be e^{{\rm i}\xi L}g_{(\xi)}(L)\,d\xi
+{\rm i}\int_0^\be g_{(\eta)}(L)e^{{\rm i}\eta L}\,d\eta.
$$
Since
$$
\int_0^\be\|e_\xi\|_{L^\be}\big\|g_{(\xi)}\big\|_{L^\be}\,d\xi<\be,
$$
it suffices to show that
$$
\lim_{s\to0}\big\|\exp({\rm i}\xi L_s)-\exp({\rm i}\xi L)\big\|=0\quad\mbox{and}\quad
\lim_{s\to0}\big\|g_{(\eta)}(L_s)-g_{(\eta)}(L)\big\|=0,\quad\xi,\,\eta>0.
$$
However, this is an immediate consequence of Corollaries \ref{ex} and Lemma \ref{sle}. $\bl$

\

\section{\bf Hilbert--Schmidt perturbations}
\setcounter{equation}{0}
\label{HS}

\

In this section we obtain estimates for $f(L)-f(M)$ in the case when $L$ and $M$ are maximal dissipative operators whose difference belongs to the Hilbert--Schmidt class $\bS_2$. We prove that all Lipschitz functions $f$ that are analytic in the upper half-plane are {\it Hilbert--Schmidt Lipschitz}, i.e.,
$$
\|f(L)-f(M)\|_{\bS_2}\le\const\|L-M\|_{\bS_2}.
$$

Put
$$
\D\df\{(x,y)\in\R^2:x=y\}.
$$

\begin{lem}
\label{diag}
Let $L$ and $M$ be maximal dissipative operators
such that $L-M\in\bS_2$. Then
$$
\iint_{\R^2}\Phi(x,y)\,d\E_L(x)(L-M)\,d\E_M(y)=0
$$
for every bounded Borel function $\Phi$ vanishing outside $\D$.
\end{lem}

\Pf By Lemma \ref{61}, for every $t\in\R$,
$$
\E_L(\{t\})L=L\E_L(\{t\})=t\E_L(\{t\})\quad\mbox{and}\quad
\E_M(\{t\})M=M\E_M(\{t\})=t\E_M(\{t\}).
$$
Besides, the sets $\big\{t\in\R:\E_L(\{t\})\not=0\big\}$ and $\big\{t\in\R:\E_M(\{t\})\not=0\big\}$
are at most countable. Hence,
\begin{align*}
\iint_{\R^2}\Phi(x,y)\,d\E_L(x)(L&-M)\,d\E_M(y)=\sum_{t\in\R}\Phi(t,t)\E_L(\{t\})(L-M)\E_M(\{t\})\\[.2cm]
&=\sum_{t\in\R}t\Phi(t,t)\big(\E_L(\{t\})\E_M(\{t\})-\E_L(\{t\})\E_M(\{t\})\big)=0.\quad\bl
\end{align*}

Let $\Li_A$ denote the set of all functions $f$ analytic
in $\C_+$ and such that \lb $f^\prime\in H^\be(\C_+)$.
Clearly, each function $f$ in $\Li_A$ extends to a function continuous on $\overline\C_+$.

\begin{thm}
\label{OLS2}
Let $f\in\Li_A\cap H^\be(\C_+)$ and let $L$ and $M$ be maximal dissipative operators
such that $L-M\in\bS_2$. Then
\bay
\label{BS2}
f(L)-f(M)=\iint_{\R^2\setminus\D}\frac{f(x)-f(y)}{x-y}\,d\E_L(x)(L-M)\,d\E_M(y),
\ey
and $\|f(L)-f(M)\|_{\bS_2}\le\|f^\prime\|_{H^\be}\|L-M\|_{\bS_2}$.
\end{thm}

\Pf In the case when $f\in \big(B_{\be1}^1(\R)\big)_+$ the result follows
from Theorem \ref{OLB} and Lemma \ref{diag}. We can take a nonnegative function $\Phi\in L^1(\R)$ such that
$\F\Phi\in C^\be(\R)$, $\supp\F\Phi\subset[-1,1]$ and $(\F\Phi)(0)=1$.
Put $\Phi_\e(x)\df\e^{-1}\Phi(\e^{-1}x)$. It is easy to see that $\Phi_\e*f\in\big(B_{\be1}^1(\R)\big)_+$
and $\|(\Phi_\e*f)^\prime\|_{H^\be}\le\|f^\prime\|_{H^\be}$ for all $\e>0$.
Moreover, $\lim_{\e\to0+}(\Phi_\e*f)(x)=f(x)$ for all $x\in\R$.
Hence,
\bey
(\Phi_\e*f)(L)-(\Phi_\e*f)(M)
=\iint_{\R^2\setminus\D}\frac{(\Phi_\e*f)(x)-(\Phi_\e*f)(y)}{x-y}\,d\E_L(x)(L-M)\,d\E_M(y)
\eey
by Theorem \ref{OLB} and Lemma \ref{diag} for all $\e>0$. It remains to observe that Lemma 3.2 in \cite{Pe6} allows us to pass to the limit as $\e\to0+$.
$\bl$

\begin{cor}
\label{OLS2C}
Let $f\in\Li_A\cap H^\be(\C_+)$ and let $L$ and $M$ be maximal dissipative operators
such that $L-M\in\bS_2$. Then
$$
f(L)-f(M)=\iint_{\R^2}\frac{f(x)-f(y)}{x-y}\,d\E_L(x)(L-M)\,d\E_M(y),
$$
and $\|f(L)-f(M)\|_{\bS_2}\le\|f^\prime\|_{H^\be}\|L-M\|_{\bS_2}$.
\end{cor}

Note that functions in $\Li_A$ do not have to be bounded and we have not defined
unbounded functions of maximal dissipative operators. However, for functions in $\Li_A$ and for maximal dissipative operators $L$ and $M$ with $L-M\in\bS_2$, we can define the operator
$f(L)-f(M)$. We cannot just put
$$
f(L)-f(M)\df\sum_{n\in\Z}\big(f_n(L)-f_n(M)\big)
$$
as we did in the previous section. Indeed, even in the scalar case we cannot write
\bay
\label{62}
f(x)-f(y)=\sum_{n\in\Z}\big(f_n(x)-f_n(y)\big)
\ey
for arbitrary $f\in\Li_A$. Formula \rf{62} can be modified in the following way.

It can be shown that for every $f\in\Li_A$, there exists a number $a$ and a sequence $\{N_j\}$ in $\Z$ such that
$$
\lim_{j\to\be}N_j=-\be\qm
f(z)-f(w)=az-aw+\lim_{j\to\be}\sum_{n=N_j}^\be(f_n(z)-f_n(w))
$$
for every $z,w\in\clos\C_+$.
This allows us to define $f(L)-f(M)$ by the formula
$$
f(L)-f(M)\df a(L-M)+\lim_{j\to\be}\sum_{n=N_j}^\be(f_n(L)-f_n(M))
$$
and the limit exists in $\bS_2$.

We do not prove this in this paper. Instead we give a different definition of  $f(L)-f(M)$ that is based on the following lemma.

\begin{lem}
There exists a sequence $\{\f_n\}_{n\ge1}$ in $H^\be(\C_+)$ such that

{\rm (i)} $\lim\limits_{n\to\be}\f_n(z)=1$ for every $z\in\C_+$,

{\rm (ii)} $\|\f_n\|_{H^\be}=1$ for every $n$,

{\rm (iii)} $({\rm i}+z)\,\f_n\in H^\be$ for every $n$,

{\rm (iv)}
$\lim\limits_{n\to\be}\|({\rm i}+z)\,\f_n^\prime(z)\|_{H^\be}=0$.
\end{lem}

\Pf
Put
$$
\f_n(z)\df\frac1{\log n}\int_1^n\frac{{\rm i}\,dt}{z+{\rm i}t}=\frac1{\log 
n}\log\frac{z+{\rm i} n}{z+{\rm i}},\quad n\ge2.
$$
Here $\log$ denotes the principal branch of logarithm.
Statements (i) and (iii) are obvious.
We have $\f_n(0)=1$ and
$$
|\f_n(z)|\le\frac1{\log n}\int_1^n\frac{dt}{|z+{\rm i}t|}\le\frac1{\log 
n}\int_1^n\frac{dt}{t}=1
$$
for all $z\in\C_+$. Hence, $\|\f_n\|_{H^\be}=1$ for every $n\ge2$.

It remains to verify (iv).
We have
$$
\big|({\rm i}+z)\,\f_n^\prime(z)\big|=\frac1{\log 
n}\cdot\left|\frac{n-1}{z+{\rm i} n}\right|
\le\frac1{\log n}
$$
for all $z\in\C_+$.
$\bl$

\begin{cor}
Let $f\in\Li_A$. Then $\f_n f\in H^\be$ for every $n$ and
$$
\lim_{n\to\be}\big\|(\f_n f)^\prime\big\|_{H^\be}=\|f^\prime\|_{H^\be}.
$$
\end{cor}

\Pf We have
$$
\|(\f_n f)^\prime\|_{H^\be}\le\|\f_n^\prime f\|_{H^\be}+\|\f_n 
f^\prime\|_{H^\be}
\le\|\f_n^\prime f\|_{H^\be}+\|f^\prime\|_{H^\be}
$$
Taking into account the fact that $|f(z)|\le\const|{\rm i}+z|$, we deduce 
from
(iv) that
$\limsup\limits_{n\to\be}\|(\f_n 
f)^\prime\|_{H^\be}\le\|f^\prime\|_{H^\be}$.
The inequality $\liminf\limits_{n\to\be}\|(\f_n 
f)^\prime\|_{H^\be}\ge\|f^\prime\|_{H^\be}$
follows from (i). $\bl$

Let $f\in\Li_A$ and let $L$ and $M$ be maximal dissipative operators
such that \lb $L-M\in\bS_2$. We can define now $f(L)-f(M)$ as follows:
\bay
\label{S2lim}
f(L)-f(M)\df\lim_{n\to\be}\big((\f_n f)(L)-(\f_n f)(M)\big)\quad \text{in the norm of}\quad\bS_2.
\ey

\begin{thm}
\label{HSLi}
Let $f\in\Li_A$. Suppose that $L$ and $M$ are maximal dissipative operators
such that $L-M\in\bS_2$. Then the limit in {\em\rf{S2lim}} exists,
$$
f(L)-f(M)=\iint_{\R^2}\frac{f(x)-f(y)}{x-y}\,d\E_L(x)(L-M)\,d\E_M(y)
$$
and
\bay
\label{nervo}
\|f(L)-f(M)\|_{\bS_2}\le\|f^\prime\|_{H^\be}\|L-M\|_{\bS_2}.
\ey
\end{thm}

\Pf We have
\begin{align*}
f(L)-f(M)&=\lim_{n\to\be}\big((\f_n f)(L)-(\f_n f)(M)\big)\nonumber\\[.2cm]
&=\lim_{n\to\be}
\iint_{\R^2}\frac{\f_n(x)f(x)-\f_n(y)f(y)}{x-y}\,d\E_L(x)(L-M)\,d\E_M(y)\nonumber\\[.2cm]
&=\iint_{\R^2}\frac{f(x)-f(y)}{x-y}\,d\E_L(x)(L-M)\,d\E_M(y),
\end{align*}
the last equality being a consequence of Lemma 3.2 in \cite{Pe6}.

This immediately implies inequality \rf{nervo}. $\bl$
 
\

\section{\bf H\"older classes and general moduli of continuity}
\setcounter{equation}{0}
\label{MC}

\

In this section we obtain estimates for $\|f(L)-f(M)\|$ for maximal dissipative operators $L$ and $M$ whose difference is bounded and for functions $f$ in the H\"older class $\big(\L_\a(\R)\big)_+$, $0<\a<1$. We show that in this case
$$
\|f(L)-f(M)\|\le\const\|L-M\|^\a,
$$
i.e., such functions $f$ are {\it operator H\"older} of order $\a$. As before we have a problem how to interpret $f(L)-f(M)$ in the case when the function $f$ is unbounded. We give the following definition
of $f(L)-f(M)$
\bay
\label{fLfM}
f(L)-f(M)\df\sum_{n=-\be}^\be\big(f_n(L)-f_n(M)\big),
\ey
where the functions $f_n$ are defined by \rf{wn}. As in the case of self-adjoint operators (see \cite{AP2}), the series on the right of converges absolutely and the definition does not depend on the choice of the functions $W_n$.

Then we proceed to the problem of estimating the operator differences $\|f(L)-f(M)\|$ for functions $f$ in the space $\big(\L_\o(\R)\big)_+$ in the case of arbitrary moduli of continuity $\o$.

\begin{thm}
\label{OHd}
There is a constant $c>0$ such that for every $\a\in(0,1)$, for arbitrary
$f\in\big(\L_\a(\R)\big)_+$, and for arbitrary maximal dissipative operators $L$ and $M$ with bounded $L-M$, the following inequality holds:
$$
\|f(L)-f(M)\|\le c\,(1-\a)^{-1}\|f\|_{\L_\a(\R)}\|L-M\|^\a,
$$
where $f(L)-f(M)$ is defined by {\em\rf{fLfM}} and the series on the right of {\em\rf{fLfM}} converges absolutely.
\end{thm}

\Pf Corollary \ref{sle} allows us to prove Theorem \ref{OHd} by using exactly 
the same arguments as in the proof of Theorem 4.1 in \cite{AP2}. $\bl$

We proceed now to the case of arbitrary moduli of continuity. For  a modulus of continuity $\o$, we define the function $\o_*$ by
\bay
\label{o*}
\o_*(x)=x\int_x^\be\frac{\o(t)}{t^2}\,dt,\quad x>0.
\ey
Then $\o_*$ is a modulus of continuity provided $\o_*(x)<\be$, $x>0$. Clearly, if $\o_*(x)<\be$ for some $x>0$, then $\o_*(x)<\be$ for all $x>0$.


For maximal dissipative operators $L$ and $M$ with bounded difference and for a function $f$ in 
$\big(\L_\o(\R)\big)_+$, we consider the operator difference $f(L)-f(M)$ and in the case of unbounded $f$, we understand by $f(L)-f(M)$ the operator
\bay
\label{LMo}
f(L)-f(M)=\sum_{n=-\be}^N\big(f_n(L)-f_n(M)\big)+\big((f-f*V_N)(L)-(f-f*V_N),(M)\big),
\ey
where the functions $V_n$ are defined by \rf{dlvp}.
As in the case of self-adjoint operators (see \cite{AP3}), the series on the right converges in the norm and the right-hand side of \rf{LMo} does not depend on the choice of the functions $W_n$.

\begin{thm}
\label{amc}
There exists a constant $c>0$ such that for every
modulus of continuity $\o$ with finite $\o_*$, every $f\in \L_\o(\R)$ and for arbitrary maximal dissipative operators $L$ and $M$ with bounded difference,
the series on the right of {\em\rf{LMo}} converges absolutely and the following inequality holds
$$
\|f(L)-f(M)\|\le c\,\|f\|_{\L_\o(\R)}\,\o_*\big(\|L-M\|\big),
$$
where the the operator $f(L)-f(M)$ is defined by {\em\rf{LMo}}.
\end{thm}

\Pf By utilizing Corollary \ref{sle}, we can prove Theorem \ref{amc}
in exactly the same way as it is done in the proof of Theorem 7.1 of \cite{AP2}. $\bl$

\begin{cor}
Let $\o$ be a modulus of continuity such that 
$$
\o_*(x)\le\const\,\o(x),\quad x>0.
$$
Then for an arbitrary function $f\in\big(\L_\o(\R)\big)_+$ and for arbitrary maximal dissipative operators $L$ and $M$ with bounded difference, the following inequality holds:
\bay
\label{n*}
\|f(L)-f(M)\|\le\const\|f\|_{\L_\o(\R)}\,\o\big(\|L-M\|\big).
\ey
\end{cor}

\

\section{\bf Higher order operator differences}
\setcounter{equation}{0}
\label{HZ}

\

In this section we establish a formula for higher operator differences in terms of multiple operator integrals. Then we obtain estimates of higher operator differences for functions of classes $\big(\L_\a(\R)\big)_+$ and
$\big(\L_{\o,m}(\R)\big)_+$.

Let $L$ and $M$ be maximal dissipative operators such that the operator $L-M$ is bounded. Put
\bay
\label{KML}
K=\frac1m(M-L).
\ey
Then the operator $L+jK$ is maximal dissipative for $0\le j\le m$. For a function $f\in\big(B_{\be1}^m(\R)\big)_+$, we consider the following finite difference
\bay
\label{konr}
\big(\D_K^m f\big)(L)\df\sum_{j=0}^m(-1)^{m-j}\left(\begin{matrix}m\\j\end{matrix}\right)f(L+jK).
\ey
By the right-hand side of \rf{konr}, we mean the following
\bay
\label{ryad}
\sum_{n\in\Z}\left(\sum_{j=0}^m(-1)^{m-j}\left(\begin{matrix}m\\j\end{matrix}\right)f_n(L+jK)\right)
\ey
where the functions $f_n$ are defined by \rf{wn}. As before, the definition does not depend on the choice of $W_n$.

The next theorem this series converges absolutely in the norm.

\begin{thm}
\label{raznk}
Let $m$ be a positive integer and let $f\in\big(B_{\be1}^m(\R)\big)_+$.
Suppose that $L$ and $M$ are maximal dissipative operators with bounded $L-M$  and let $K$ be the operator defined by {\em\rf{KML}}.
Then the series {\em\rf{ryad}} converges absolutely in the norm and
$$
\big(\D_K^m f\big)(L)=m!\int\cdots\int(\dg^mf)(x_1,\cdots,x_{m+1})
\,d\E_1(x_1)K\,d\E_2(x_2)K\cdots K\,d\E_{m+1}(x_{m+1}),
$$
where $\E_j$ is the semi-spectral measure of $L+jK$.
\end{thm}

Theorem \ref{raznk} implies the following result:

\begin{thm}
\label{Kn}
Let $m$ be a positive integer. There exists $c>0$ such that for arbitrary $f$, $L$, and $M$ satisfying the hypotheses of Theorem {\em\ref{raznk}} the following inequality holds:
$$
\big\|\big(\D_K^m f\big)(L)\big\|\le c\,\|f\|_{B_{\be1}^m(\R)}\|K\|^m.
$$
\end{thm}

\Pf The result follows immediately from Theorem \ref{raznk} and from Theorem 5.5 of \cite{Pe5}. $\bl$

\medskip

To avoid complicated notation, we prove Theorem \ref{raznk} in the case $m=2$. The proof can easily be adjusted for an arbitrary positive integer $m$.

\begin{lem}
\label{in2}
Let $L$ be a maximal dissipative operator and let $K$ be a bounded operator
such that $L+2K$ is a dissipative operator. Then
\begin{align*}
\exp\big({\rm i}x(L&+2K))-2\exp({\rm i}x(L+K))+\exp({\rm i}xL\big)\\[.2cm]
=&-2\int_0^x\left(\int_0^\xi\exp({\rm i}(x-\xi)L)K\exp({\rm i}
(\xi-\eta)(L+K))K\exp({\rm i}\eta(L+2K))\,d\eta\right)\,d\xi
\end{align*}
for every $x>0$.
\end{lem}

\Pf By
Corollary \ref{in}, for $x>0$, we have
$$
\exp\big({\rm i}x(L+2K)\big)-\exp({\rm i}xL)=
2{\rm i}\int_0^x\exp\big({\rm i}(x-\xi)L\big)K\exp\big({\rm i}\xi (L+2K)\big)\,d\xi
$$
and
$$
\exp\big({\rm i}x(L+K)\big)-\exp({\rm i}xL)=
{\rm i}\int_0^x\exp\big({\rm i}(x-\xi)L\big)K\exp\big({\rm i}\xi (L+K)\big)\,d\xi.
$$
It also follows from Corollary \ref{in} that for $\xi>0$,
$$
\exp\big({\rm i}\xi(L+2K)\big)-\exp\big({\rm i}\xi(L+K)\big)=
{\rm i}\int_0^\xi\exp\big({\rm i}(\xi-\eta)(L+K)\big)K\exp\big({\rm i}\eta (L+2K)\big)\,d\eta.
$$
These equalities imply the result. $\bl$

\begin{lem}
\label{H1+}
Let $f$ be a function in $H^1(\C_+)$ such that
$$
\int_0^\be(1+|\xi|^2)\big|(\F f)(\xi)\big|\,d\xi<\be.
$$
Then under the hypotheses of Lemma {\em\ref{in2}},
\begin{align*}
&f(L+2K)-2f(L+K)+f(L)\\[.2cm]
=&-\frac1\pi\iiint\limits_{\R_+^3}(\F f)(\xi+\eta+\vk)\exp(\ri\xi L)K
\exp\big(\ri\eta (L+K)\big)K\exp\big(\ri\vk (L+2K)\big)\,d\xi\,d\eta\,d\vk.
\end{align*}
\end{lem}

\Pf The result follows from Lemma \ref{in2} in the same way as Lemma \ref{H1}
has been deduced from Corollary \ref{in}. Indeed, we have
\begin{align*}
f(L+2K)-&2f(L+K)+f(L)\\[.2cm]
=&\frac1{2\pi}\int_0^\be(\F f)(x)\Big(\exp({\rm i}x(L+2K))-2\exp({\rm i}x(L+K))+\exp({\rm i}xL)\Big)\,dx.
\end{align*}
It remains to apply Lemma \ref{in2} and to change variables. $\bl$

Now we are ready to prove Theorem \ref{raznk} for $m=2$.

\begin{thm}
\label{OLB2}
Let $f\in \big(B_{\be1}^2(\R)\big)_+$, let $L$ be a maximal dissipative operator
and let $K$ be a bounded operator
such that $L+2K$ is dissipative. Then
\bey
f(L+2K)-2f(L+K)+f(L)=2\iiint(\dg^2f)(x,y,z)d\E_L(x)Kd\E_{L+K}(y)Kd\E_{L+2K}(z).
\eey
\end{thm}

\Pf As in the proof of Theorem \ref{OLB}, it suffices to consider the case when $f\in L^\be(\R)$
and the support of its Fourier transform $\F f$ is a compact subset of $(0,\be)$.
Then we have (see \cite{Pe5})
\begin{align}
\label{tensor}
(\dg^2f)(x,y,z)=&
-\iint\limits_{\R_+\times\R_+}f_{(\eta+\vk)}(x)e^{{\rm i}\eta y}e^{{\rm i}\vk z}
\,d\eta\,d\vk\nonumber\\[.2cm]
&-\iint\limits_{\R_+\times\R_+}e^{{\rm i}\xi x}f_{(\xi+\vk)}(y)e^{{\rm i}\vk z}
\,d\xi\,d\vk
-\iint\limits_{\R_+\times\R_+}e^{{\rm i}\xi x}e^{{\rm i}\eta y}f_{(\xi+\eta)}(z)
\,d\xi\,d\eta.
\end{align}
The rest of the proof ids the same as in the proof of Lemma \ref{OLe}. $\bl$

We proceed now to estimation of operator finite differences.

\begin{thm}
\label{alfam}
Let $\a$ be a positive integer and let $m$ be an integer such that $\a<m$. There exists $c>0$ such that for arbitrary dissipative operators $L$ and $M$ with bounded $L-M$ and for every $f\in\big(\L_\a(\R)\big)_+$ the following inequality holds:
$$
\big\|\big(\D_K^m f\big)(L)\big\|\le c\,\|f\|_{\L_\a(\R)}\|K\|^\a,
$$
where $K$ is defined by {\em\rf{KML}}.
\end{thm}

\Pf It follows from Theorem \ref{raznk} that
\bay
\label{2nm}
\big\|\big(\D_K^m f_n\big)(L)\big\|\le\const2^{nm}\|f_n\|_{L^\be}\|K\|^m,
\ey
where $f_n$ is defined by \rf{wn}.

The rest of the proof is the same as in the proof of Theorem 4.2 of \cite{AP2}. $\bl$

The following theorem yields estimates of $\big\|\big(\D_K^m f_n\big)(L)\big\|$ for functions $f$ in $\big(\L_{\o,m}(\R)\big)_+$, where $\o$ is a nondecreasing function satisfying \rf{on}. With such a function $\o$ we associate the auxiliary function 
$\o_{*,m}$ defined by
$$
\o_{*,m}(x)=x^m\int_x^\be\frac{\o(t)}{t^{m+1}}\,dt=\int_1^\be\frac{\o(sx)}{s^{m+1}}\,dx.
$$

\begin{thm}
\label{Lom}
Let $m$ be a positive integer.  Then there is a positive number $c$ such that for an arbitrary
nondecreasing function $\o$ on $(0,\be)$ satisfying {\em\rf{on}} with finite $\o_{*,m}$, for an arbitrary
function $f$ in $\big(\L_{\o,m}(\R)\big)_+$, and arbitrary maximal dissipative operators $L$ and 
$M$ with bounded $L-M$ the following inequality holds:
$$
\left\|\big(\D_K^mf\big)(L)\right\|\le c\,\|f\|_{\L_{\o,m}(\R)}\,\o_{*,m}\big(\|K\|\big),
$$
where $K$ is defined by {\em\rf{KML}}.
\end{thm}

\Pf If we apply inequality \rf{2nm}, we can proceed in the same way as in the proof of Theorem 7.1 of \cite{AP4}. $\bl$

\begin{cor}
Suppose that in the hypotheses of Theorem {\em\ref{Lom}} the function $\o$ satisfies the condition
$$
\o_{*,m}(x)\le\const\o(x),\quad x>0.
$$
Then
$$
\left\|\big(\D_K^mf\big)(L)\right\|\le\const\|f\|_{\L_{\o,m}(\R)}\,\o\big(\|K\|\big).
$$
\end{cor}


\

\section{\bf Higher operator derivatives}
\setcounter{equation}{0}
\label{dif}

\

In this section we show that under the assumption $f\in\big(B_{\be1}^m(\R)\big)_+$, the function
\bay
\label{fotLt}
t\mapsto f(L_t)
\ey
has $m$th derivative and we obtain a formula the $m$th derivative in terms of multiple operator integrals. 

Here $L$ and $M$ are maximal dissipative operators such that the operator $L-M$ is bounded and the family $L_t$ is defined by \rf{Lt}.

To be more precise, we note that functions in $\big(B_{\be1}^m(\R)\big)_+$ do not have to be bounded and in the case $m>1$ under the assumption $f\in\big(B_{\be1}^m(\R)\big)_+$, the function \rf{fotLt} does not have to be differentiable. However, we show that for each $n\in\Z$, the function
$$
t\mapsto f_n(L_t)
$$
has derivatives in the norm up to order $m$ and the series
\bay
\label{sumfn}
\sum_{n\in\Z}\frac{d^m}{dt^m}f_n(L_t)
\ey
converges absolutely in the norm and be the $m$th derivative of the function \rf{fotLt} we understand the sum of the series \rf{sumfn}.

\begin{thm}
\label{proiz}
Let $m$ be a positive integer and let $f\in\big(B_{\be1}^m(\R)\big)_+$.
Suppose that $L$ and $M$ are maximal dissipative operators such that
$L-M$ is bounded. Then the function {\em\rf{fotLt}} has $m$th derivative and
$$
\frac{d^m}{dt^m}f(L_t)\Big|_{t=s}\!=m!\!
\int\!\!\cdots\!\!\int\big(\dg^mf)(x_1,x_2,\cdots,x_{m+1})
\,d\E_s(x_1)K\,d\E_s(x_2)K\cdots K\,d\E_s(x_{m+1}),
$$
where $\E_s$ is the semi-spectral measure of $L_s$.
\end{thm}

\Pf Clearly, it suffices to prove the result for $s=0$.
It follows from \rf{Boke} that 
\begin{align*}
\sum_{n\in\Z}&\left\|
\int\cdots\int\big(\dg^mf_n)(x_1,x_2,\cdots,x_{m+1})
\,d\E(x_1)K\,d\E(x_2)K\cdots K\,d\E(x_{m+1})\right\|\\[.2cm]
&\le\sum_{n\in\Z}2^{nm}\|f_n\|_{L^\be}\|K\|^m
\le\const\|f\|_{B_{\be1}^m(\R)}\|K\|^m,
\end{align*}
where $\E$ is the semi-spectral measure of $L$.
Hence, it is sufficient to prove that
$$
\frac{d^m}{dt^m}f_n(L_t)\Big|_{t=0}\!=m!\!
\int\!\!\cdots\!\!\int\big(\dg^mf_n)(x_1,x_2,\cdots,x_{m+1})
\,d\E(x_1)K\,d\E(x_2)K\cdots K\,d\E(x_{m+1}).
$$
To simplify the notation, we prove this identity for $m=2$. In the general case the proof is the same. 

Put $g=f_n$. We need the following identities:
\begin{align*}
&\frac1t\left(\iint\big(\dg f_n\big)(x,y)\,d\E_t(x)K\,d\E_t(y)-
\iint\big(\dg f_n\big)(x,y)\,d\E_t(x)K\,d\E(y)\right)\\
=&\iiint\big(\dg^2f_n\big)(x,y,z)\,d\E_t(x)K\,d\E_t(y)K\,d\E(z)
\end{align*}
and
\begin{align*}
&\frac1t\left(\iint\big(\dg f_n\big)(x,y)\,d\E_t(x)K\,d\E(y)
-\iint\big(\dg f_n\big)(x,y)\,d\E_A(x)K\,d\E(y)\right)\\
=&\iiint\big(\dg^2f_n\big)(x,y,z)\,d\E_t(x)K\,d\E(y)K\,d\E_t(z).
\end{align*}

The proof of these identities is similar to the proof of Theorem \ref{OLB2}.

It follows that
\begin{align*}
\frac1t\left(\frac{d}{ds}f_n(L_s)\Big|_{s=t}-\frac{d}{ds}f(L_s)\Big|_{s=0}\right)
&=\iiint\big(\dg^2f_n\big)(x,y,z)\,d\E_t(x)K\,d\E_t(y)K\,d\E(z)\\
&+\iiint\big(\dg^2f_n\big)(x,y,z)\,d\E_t(x)K\,d\E(y)K\,d\E(z).
\end{align*}
It remains to observe that
\begin{align*}
\lim_{t\to0}&
\iiint\big(\dg^2f_n\big)(x,y,z)\,d\E_t(x)K\,d\E_t(y)K\,d\E(z)\\[.2cm]
=&\iiint\big(\dg^2f_n\big)(x,y,z)\,d\E(x)K\,d\E(y)K\,d\E(z)
\end{align*}
and
\begin{align*}
\lim_{t\to0}&
\iiint\big(\dg^2f_n\big)(x,y,z)\,d\E_t(x)K\,d\E(y)K\,d\E(z)\\[.2cm]
=&\iiint\big(\dg^2f_n\big)(x,y,z)\,d\E(x)K\,d\E(y)K\,d\E(z).
\end{align*}
This can be proved in the same way as in the proof of Theorem \ref{der}
if we apply \rf{tensor}. $\bl$

\

\section{\bf Estimates in Schatten--von Neumann classes}
\setcounter{equation}{0}
\label{SvN}

\

In this section we consider the case when $f\in\big(\L_\a(\R)\big)_+$, $\a>0$, and $L$ and $M$ are maximal dissipative operators such that $L-M$ belongs to the Schatten--von Neumann class $\bS_p$. We are going to obtain results that are similar to the results of \cite{AP3} for self-adjoint operators. 

We need the following fact:

\begin{lem}
\label{spm} 
Let $f\in\big(B_{\be1}^m(\R)\big)_+$ and let $p\ge m$. Then there exists a positive number $c$ such that for arbitrary maximal dissipative operators 
$L$ and $M$  such that $L-M$ belongs to the Schatten--von Neumann class $\bS_p$ and for every $H^\be(\C_+)$ function $f$ whose Fourier transform is supported on $[0,\s]$, the following inequality holds:
$$
\big\|\big(\D_K^m f\big)(L)\big\|_{\bS_{p/m}}
\le c\,\s^m\|f\|_{L^\be}\|K\|^m_{\bS_p},
$$
where $K$ is defined by {\em\rf{KML}}.
\end{lem}

\Pf The result follows easily from Theorem \ref{raznk} and \rf{Boke}. $\bl$

\begin{thm}
\label{Spa}
Let $\a>0$ and let $m$ be an integer such that $m-1\le\a<m$. Suppose that $p>m$. Then there exists a positive number $c$ such that for arbitrary maximal dissipative operators $L$ and $M$ with $L-M\in\bS_p$ and for every $f\in\big(\L_\a(\R)\big)_+$, the operator
$\big(\D_K^m f\big)(L)$ belongs to $\bS_{p/\a}$ and 
$$
\big\|\big(\D_K^m f\big)(L)\big\|_{\bS_{p/\a}}\le
c\,\|f\|_{\L_\a(\R)}\|K\|_{\bS_p}^m.
$$
\end{thm}

In the case $p=m$ we can obtain a weaker conclusion that $\big(\D_K^m f\big)(L)$ belongs to the ideal $\bS_{\frac m\a,\be}$. Recall that the ideal $\bS_{q,\be}$ consists of operators $T$ on Hilbert space whose singular values $s_j(T)$ satisfy the condition:
$$
s_j(T)\le\const(1+j)^{-1/q}.
$$

\begin{thm}
\label{Spbe}
Let $\a>0$ and let $m$ be an integer such that $m-1\le\a<m$. Then there exists a positive number $c$ such that for arbitrary maximal dissipative operators $L$ and $M$ with $L-M\in\bS_m$ and for every $f\in\big(\L_\a(\R)\big)_+$, the operator
$\big(\D_K^m f\big)(L)$ belongs to $\bS_{\frac m\a,\be}$ and 
$$
\big\|\big(\D_K^m f\big)(L)\big\|_{\bS_{\frac m\a,\be}}\le
c\,\|f\|_{\L_\a(\R)}\|K\|_{\bS_m}^m.
$$
\end{thm}

To prove Theorems \ref{Spa} and \ref{Spbe}, we can use Lemma \ref{spm} and
then proceed along the same line as it was done in \S\,5 of \cite{AP3}.

\medskip

{\bf Remark.} Note that in the case $m=1$ the assumptions of Theorem \ref{Spbe} do not guarantee that $\big(\D_K^m f\big)(L)\in\bS_{m/\a}$. Indeed, in \S\,9 of \cite{AP3} it was shown that the corresponding fact for self-adjoint operators does not hold. The example given in \cite{AP3} is based on the Schatten--von Neumann criterion for Hankel operators, see \cite{Pe0} and \cite{Pe4}. On the other hand, if $m$ is an integer greater than 1, we do not know whether the assumptions of Theorem \ref{Spbe} imply that $\big(\D_K^m f\big)(L)\in\bS_{m/\a}$.

The following result shows that if we impose a slightly stronger assumption on $f$, we can obtain the conclusion that $\big(\D_K^m f\big)(L)\in\bS_{m/\a}$.

\begin{thm}
\label{SpBes}
Let $\a>0$ and let $m$ be an integer such that $m-1\le\a\le m$. Then there exists a positive number $c$ such that for arbitrary maximal dissipative operators $L$ and $M$ with $L-M\in\bS_m$ and for every $f\in\big(B_{\be1}^\a(\R)\big)_+$, the operator
$\big(\D_K^m f\big)(L)$ belongs to $\bS_{\frac m\a,\be}$ and 
$$
\big\|\big(\D_K^m f\big)(L)\big\|_{\bS_{m/\a}}\le
c\,\|f\|_{B_{\be1}^\a(\R)}\|K\|_{\bS_m}^m.
$$
\end{thm}

We can improve Theorem \ref{Spa} in the following way.

\begin{thm}
\label{chsSp}
Let $\a>0$, $m-1\le\a<m$, and $m<p<\be$. Then there exists a positive number $c$ such that for every $f\in\big(\L_\a(\R)\big)_+$, every $l\in\Z_+$, and for arbitrary maximal dissipative operators $L$ and $M$ with bounded $L-M$, the following inequality holds:
$$
\sum_{j=0}^l\left(s_j\left(\big|\big(\D_K^mf\big)(L)\big|^{1/\a}\right)\right)^p
\le c\,\|f\|_{\L_\a(\R)}^{p/\a}\sum_{j=0}^l\big(s_j(K)\big)^p.
$$
\end{thm}

Again, if we use Lemma \ref{spm}, we can prove Theorems \ref{SpBes} and
\ref{chsSp} in the same way as the proofs of the corresponding facts for self-adjoint operators given in \S\,5 of \cite{AP3}.

Note also that as in the case of self-adjoint operators (see \cite{AP3}), we can obtain more general results for ideals of operators on Hilbert space with upper Boyd index less than 1.

\

\section{\bf Commutators and quasicommutators}
\setcounter{equation}{0}
\label{cq}

\

In this section we estimate quasicommutators $f(L)R-Rf(M)$ in terms of $LR-RM$, where $L$ and $M$ are maximal dissipative operators and $R$ is a bounded operator. In the special case $R=I$ we arrive at the problem of estimating $f(L)-f(M)$ in terms of $L-M$; this problem was discussed in previous sections. On the other hand, in the case $L=M$ we get the problem of estimating commutators $f(L)R-Rf(L)$.

Let us explain what we mean by the boundedness of $LR-RM$
for not necessarily bounded normal operators $L$ and $M$.

We say that the {\it operator $LR-RM$ is bounded} if
$R(\cd_M)\subset \cd_{L}$ and 
$$
\|LRu-RMu\|\le\const\|u\|\quad\mbox{for every}\quad u\in \cd_M.
$$
Then there exists a unique bounded operator $K$ such that
$Ku=LRu-RMu$ for all $u\in \cd_{M}$. In this case we write $K=LR-RM$. Thus $N_1R-RN_2$ is bounded if and only if
$$
\big|(Ru,L^*v)-(Mu,R^*v)\big|\le \const\|u\|\cdot\|v\|
$$
for every $u\in \cd_{M}$ and $v\in \cd_{L^*}$ (recall that $-L^*$ and $-M^*$ are also maximal dissipative operators, see \S\,\ref{dis}). It is easy to see that $LR-RM$ is bounded if and only if $M^*R^*-R^*L^*$ is bounded,
and $(LR-RM)^*=-(M^*R^*-R^*L^*)$.
In particular, we write $LR=RM$ if $R(\cd_M)\subset \cd_L$ and $LRu=RMu$ for every $u\in \cd_M$.
We say that $\|LR-RM\|=\be$ if $LR-RM$ is not a bounded operator.

\begin{thm}
\label{sigma}
Let $f$ be a function in $H^\be(\C_+)$ such that $\supp\F f\subset[0,\s]$.
Suppose that $L$ and $M$ are maximal dissipative operators and $R$ is a bounded operator such that the operator $LR-RM$ is bounded. Then
\bay
\label{qsig}
f(L)R-Rf(M)=\iint\frac{f(x)-f(y)}{x-y}\,d\E_L(x)(LR-RM)\,d\E_M(y)
\ey
and
$$
\big\|f(L)R-Rf(M)\big\|\le8\s\|LR-RM\|.
$$
\end{thm}

The proof of Theorem \ref{sigma} is similar to the proof of the corresponding result for $f(L)-f(M)$, see \S\,\ref{OLD}.

For unbounded functions $f$, we have not defined $f(L)$ and $f(M)$. However, we show in this section that under certain natural assumptions on $f$, it is possible to define the quasicommutator $f(L)R-Rf(M)$ for unbounded $f$ by the formula
\bay
\label{cryad}
\sum_{n\in\Z}\big(f_n(L)R-Rf_n(M)\big).
\ey

\begin{thm}
\label{qcom}
Let $f\in\big(B_{\be1}^m(\R)\big)_+$. Suppose that $L$ and $M$ are maximal dissipative operators and $R$ is a bounded operator such that the operator 
$LR-RM$ is bounded. Then the series {\em\rf{cryad}} converges absolutely in the norm, formula {\em\rf{qsig}} holds and
$$
\big\|f(L)R-Rf(M)\big\|\le\const\|f\|_{B_{\be1}^m(\R)}\|LR-RM\|.
$$
\end{thm}

\Pf The result follows immediately from Theorem \ref{sigma} and the definition of the Besov space $\big(B_{\be1}^m(\R)\big)_+$. $\bl$

\begin{thm}
\label{comH}
Let $0<\a<1$ and let $f\in\big(\L_\a(\R)\big)_+$. Suppose that $L$ and $M$ are maximal dissipative operator and $R$ is a bounded operator such that the operator 
$LR-RM$ is bounded. Then the series {\em\rf{cryad}} converges absolutely in the norm and
$$
\big\|f(L)R-Rf(M)\big\|\le\const\|f\|_{\L_\a(\R)}\|LR-RM\|^\a\|R\|^{1-\a}.
$$
\end{thm}

We proceed now to the case of functions in the space $\big(\L_\o(\R)\big)_+$, where $\o$ is an arbitrary modulus of continuity. For $f\in\big(\L_\o(\R)\big)_+$, we define $f(L)R-Rf(M)$ by
\begin{align}
\label{flrm}
f(L)R-Rf(M)&=\sum_{n=-\be}^N\big(f_n(L)R-Rf_n(M)\big)\nonumber\\[.2cm]
+&\big((f-f*V_N)(L)R-R(f-f*V_N)(M)\big).
\end{align}

\begin{thm}
\label{modne}
Let $\o$ be a modulus of continuity such that the function $\o_*$ defined by {\em\rf{o*}} takes finite values.
Suppose that $L$ and $M$ are maximal dissipative operator and $R$ is a bounded operator such that the operator $LR-RM$ is bounded. Then for arbitrary
$f\in\big(\L_\o(\R)\big)_+$, the series in {\em\rf{flrm}} converges absolutely in the norm and the following inequality holds:
$$
\|f(L)R-Rf(M)\|\le\const\,\|f\|_{\L_\o(\R)}
\|R\|\,\o_*
\left(\frac{\|LR-RM\|}{\|R\|}\right).
$$
\end{thm}

The proof is the same as the proof of Theorem \ref{amc}.

In \S\,\ref{SvN} we have obtained estimates for $f(L)-f(M)$ when $L-M$ belongs to the Schatten--von Neumann class $\bS_p$ (this corresponds to the case $m=1$ in \S\,\ref{SvN}). We can obtain analogs of all those results for quasicomutators. We state here an analog of Theorem \ref{Spa}.

\begin{thm}
\label{comSp}
Let $p>1$, $0<\a<1$, and let $f\in\big(\L_\a(\R)\big)_+$. Suppose that $L$ and $M$ are maximal dissipative operator and $R$ is a bounded operator such that the operator $LR-RM\in\bS_p$. Then $f(L)R-Rf(M)\in\bS_{p/\a}$ and
$$
\big\|f(L)R-Rf(M)\big\|_{\bS_{p/\a}}\le\const\|f\|_{\L_\a(\R)}\|LR-RM\|_{\bS_p}^\a\|R\|^{1-\a}.
$$
\end{thm}

The proof of Theorem \ref{comSp} is almost the same as that of Theorem \ref{Spa} in the case $m=1$.

\

\

\noindent
\begin{tabular}{p{9cm}p{15cm}}
A.B. Aleksandrov & V.V. Peller \\
St-Petersburg Branch & Department of Mathematics \\
Steklov Institute of Mathematics  & Michigan State University \\
Fontanka 27, 191023 St-Petersburg & East Lansing, Michigan 48824\\
Russia&USA
\end{tabular}

\end{document}